\documentclass[10pt,a4paper]{article}
\usepackage[utf8]{inputenc}
\usepackage[english]{babel}
\usepackage[T1]{fontenc}
\usepackage{amsmath}
\usepackage{amsfonts}
\usepackage{amsthm}
\usepackage{amsmath}
\usepackage{amssymb}
\usepackage{mathrsfs}
\usepackage{mathtools}
\usepackage[toc,page]{appendix}
\usepackage{dsfont} 
\usepackage{enumerate}
\usepackage[round]{natbib}
\usepackage[dvipsnames]{xcolor} 
\usepackage{graphicx}
\usepackage{amsmath}
\usepackage{float}
\usepackage{makecell}
\usepackage{times}
\usepackage{fancyhdr}

\usepackage{algorithm}
\usepackage{algcompatible}
\usepackage{caption}
\usepackage{subcaption}
\usepackage{tikz}
\tikzstyle{int}=[draw, line width = 0.5mm, minimum size=5em]
\usetikzlibrary{positioning,automata,arrows}
\usetikzlibrary{matrix}
\usepackage[normalem]{ulem}

\author{Henri Mermoz KOUYE$^{(1)}$ and Gildas MAZO$^{(1)}$\\
 \footnotesize{$^{(1)}$Univ.  Paris-Saclay, INRAE, MaIAGE, 78350, Jouy-en-Josas, France}}

\title{Regularizing nested Monte Carlo Sobol' index estimators to
  balance the trade-off between explorations and repetitions in
  global sensitivity analysis of stochastic models}

\newtheorem{assumption}{Assumption}
\newtheorem{proposition}{Proposition}

\newtheorem{corollary}{Corollary}

\newtheorem{lemmaa}{Lemma}[section]
\newtheorem{theorem}{Theorem}

\newtheorem{remark}{Remark}

\newcommand\del{\bgroup\markoverwith
  {\textcolor{red}{\rule[2pt]{2pt}{0.5pt}}}\ULon}

\DeclareMathOperator{\expec}{E}
\DeclareMathOperator{\var}{Var}

\newcommand{\V}[1]{\var\left(#1\right)}

\newcommand{\Vc}[2]{\var\left[#1\mid #2\right]}
\newcommand{\E}[1]{\mathbb{E}\left(#1\right)}
\newcommand{\Esmall}[1]{\mathbb{E}(#1)}
\newcommand{\Ec}[2]{\mathbb{E}\left[#1\mid #2\right]}

\newcommand{\X}{\mathbf{X}}

\usepackage{mwe}

\newlength{\tempdima}

\newlength{\tempwidth}
\newcommand{\rowname}[1]
{\rotatebox{90}{\makebox[\tempdima][c]{\textbf{#1}}}}
\newcommand{\columnname}[1]
{\makebox[\tempwidth][c]{\textbf{#1}}}

\begin{document}
 \maketitle

\begin{abstract}
  Sobol' sensitivity index estimators for stochastic models are
  functions of nested Monte Carlo estimators, which are estimators
  built from two nested Monte Carlo loops.  The outer loop explores
  the input space and, for each of the explorations, the inner loop
  repeats model runs to estimate conditional expectations.  Although
  the optimal allocation between explorations and repetitions of one's
  computational budget is well-known for nested Monte Carlo
  estimators, it is less clear how to deal with functions of nested
  Monte Carlo estimators, especially when those functions have
  unbounded Hessian matrices, as it is the case for Sobol' index
  estimators.  To address this problem, a regularization method is
  introduced to bound the mean squared error of functions of nested
  Monte Carlo estimators. Based on a heuristic, an
  allocation strategy that seeks to minimize a bias-variance trade-off
  is proposed. The method is applied to Sobol' index estimators for
  stochastic models.  A practical algorithm that adapts to the level
  of intrinsic randomness in the models is given and illustrated on
  numerical experiments.
\end{abstract}


\section{Introduction}
(Global) sensitivity analyses of stochastic models may be
challenging.  Indeed, stochastic models
include two sources of uncertainty: the parameters' uncertainty and
the intrinsic randomness of the stochastic model. The latter can be
seen as a hidden additional random input, which may challenge the
definition of meaningful sensitivity indices and their efficient
estimation.

Many methods have been introduced to deal with stochastic models. A
first approach proposed by \cite{Hart} considers random Sobol'-Hoeffding
decompositions of stochastic model outputs and defines
sensitivity indices of such models as expectations of the resulting
random Sobol' indices. A second approach, widely used in
applications~\citep{courcoul}, 
 focuses on deterministic
quantities of interest (QoIs) such as conditional expectations or
conditional variances.  By conditioning with
respect to the uncertain parameters, the aim is to smooth the
intrinsic randomness out and hence to deal with 
 deterministic functions of the uncertain parameters only, so
that sensitivity analysis (SA) methods for deterministic models can be applied. A third
approach includes recently developed
methods~\citep{fort2020global,gamboa2021sensitivity,daveiga} that see
stochastic models as deterministic functions with values in
probability distribution spaces. Various sensitivity indices are
defined on such spaces in order to measure the contributions of the
uncertain parameters. A fourth approach rests on
meta-models~\citep{viet,etore,jimenez,FORT2,zhu,janon,ivan}.
As a remark, let us note that that estimation of
expectations of functions of conditional moments can be performed with
specific methods, such as semi-parametric~\citep{daveiga2,daveiga3}
or multilevel Monte Carlo~\citep{mycek} methods.
Estimation in high-dimensional deterministic models is considered in~\cite{castro}.

In the first three approaches described
above,
estimation relies on two nested Monte Carlo loops. Indeed, not only
should the model be evaluated at many points in the input space (it is
said that the input space is explored), but also the model should be
run several times (it is said that the model is ``repeated'') at each
of those explorations to estimate conditional expectations. In the
first approach, these repetitions are performed when approximating the
expectations of the random Sobol' indices. In the second and third
approaches, model outputs are repeated when estimating the QoIs and
the probability distributions, respectively.
The larger the number of explorations and the number of repetitions,
the better the sensitivity index estimates.  However, in
practice, models may be complex and computationally
expensive. Therefore, getting a large number of runs may be
impossible, and, therefore, finding a compromise between  explorations
and repetitions under the constraint of a computational
cost or precision of estimates is important.  For instance,
\cite{mazo}  proposed a choice of the number
of explorations and the number of repetitions based on the
minimization of some bound of the so-called mean ranking error
of the estimators. This error measures the gap between the ranks of
the theoretical indices and those of the estimators.  However, a small
mean ranking error does not necessarily imply that estimations are close to their
theoretical values.

A nested Monte Carlo estimator is an estimator built from
two nested Monte Carlo loops for estimating an expectation of a
function of a conditional expectation. The compromise between the
sizes of the outer and inner loops, coinciding, respectively, with the
number of explorations and repetitions in SA, is a
bias-variance trade-off. This question is
sometimes refered to as the question of optimal allocation. Typically, the mean-squared error (MSE) is
approximated and minimized to find the ``best'' asymptotic
trade-off~\citep{rainforth,giles2019multilevel,giorgi_weak_2020,juneja}.
More precisely, if $T$ denotes the number of model runs  available to
the practitioner, and
$n=T^{1-\eta}$ and $m=T^{\eta}$, $\eta\in [0,1)$, denote the sizes of the outer and
inner loops, respectively, then it is desired to know which value of
$\eta$ would lead to the optimal convergence rate of the MSE.

Sobol' index estimators are \emph{functions} of nested Monte Carlo
estimators. It is, therefore,  natural to
ask whether the MSE of functions of nested Monte Carlo
estimators can be bounded, too.  The answer will depend on the
regularity of the function in question. For instance, if the function
is twice continuously differentiable with a bounded Hessian matrix,
then standard Taylor expansion techniques can be used. However, if the
Hessian is unbounded, as it is the case of Sobol' index
estimators, then the question is more challenging and remains open, to
the best of our knowledge. 

We note that the question of optimal
allocation of  available computational resources is shared more
globally by research fields interested in stochastic
simulators~\citep{chen1,CHEN2017575,binois1,binois2}.

In this paper, we introduce a regularized version of the nested Monte
Carlo Sobol' index estimator and bound its
MSE. Based on a heuristic, the bound is simplified and a
minimization is carried out to balance the numbers of explorations and
repetitions. A practical algorithm that
adapts to the model intrinsic randomness is proposed and illustrated
on numerical experiments.

This paper is organized as follows. Sobol' indices for stochastic
models and their estimation based on nested Monte Carlo methods are
presented in Section~\ref{sec:context}. The regularization method, as
well as the challenges of bounding
the MSE of functions of nested Monte Carlo estimators, are
introduced in Section~\ref{sec:setting}. Section
\ref{sec:application-to-sa} is an application to
 Sobol' index estimators.  Practical algorithms
are presented and illustrated in Section \ref{algo-illus-section}. A
conclusion closes the paper.


\section{Sensitivity analysis of stochastic models} 
\label{sec:context}

A stochastic model with
inputs $\X=(X_1,\dots,X_p)\in\mathbb{R}^p$ and output $Y\in \mathbb{R}$ is modeled as a function $f$
of $\X$ and some collection of random variables, denoted by $Z$,  independent of $\X$
such that
\begin{equation}
Y=f\left(\X,Z\right).\label{stoch-model}
\end{equation}
The stochasticity of the model originates from $Z$
since the output of the model evaluated at an input $\X=x$ is a random
variable $f(x,Z)$. The distribution of $Z$ is generally unknown. 

In the context of SA, a way to deal with stochastic
models consists in carrying out SA for deterministic
models given by deterministic QoIs. This allows to
switch from a stochastic model to some deterministic models for which
many SA methods are studied in the literature.

We consider QoIs of the form
\begin{equation}
Q(\X)=\Ec{\varphi(\X,Z)}{\X}, \label{qoi-form}
\end{equation}
where $\varphi(\X,Z)$ is a function of $\X$ and $Z$.  For instance, if
$\varphi=f$ then $Q(\X)$ is the conditional expectation of the model and if
$\varphi(\X,Z) = \left(f(\X,Z) - \Ec{f(\X,Z)}{\X}\right)^2$ then
$Q(\X)$ is the conditional variance, two common choices in practice.

If $u$ is a subset of $\{1,\dots,p\}$, denote by $\X_u$ the group of
inputs $\{X_i,i \in u\}$ and $\X_{\sim u}$ the group of inputs
$\{X_i,i \not\in u\}$.  The Sobol' and total indices of the input
vector $\X_u$ with respect to the function $Q$ are defined as

\begin{align}
S_u  &=\frac{\V{\Ec{Q(\X)}{\X_u}}}{\V{Q(\X)}} \\
T_u  &=1-\frac{\V{\Ec{Q(\X)}{\X_{\sim u}}}}{\V{Q(\X)}}=1-S_{\sim u}.
\end{align}

The Sobol' index $S_u$ (and hence $T_u$) can be expressed in
terms of a function $g$ linking the components of some parameter
vector. Let $\widetilde{\X}$ be an independent copy of $\X$,
independent of $Z$.  Denote by $\widetilde{\X}_{\sim u}$ the subvector
of $\widetilde\X$ whose components are those of $\widetilde\X$ not
indexed by $u$. (For instance, if $p=4$ and $u=\{1,4\}$ then
$\widetilde{\X}_{\sim u}=(\widetilde X_2,\widetilde X_3)$.)
If $\theta=(\theta_1,\theta_2,\theta_3)$ with
$\theta_1 = \Esmall{Q(\X)^2}$, $\theta_2 = \Esmall{Q(\X)}$ and
$\theta_3=\Esmall{\Ec{Q(\X)}{\X_u}^2} =
\Esmall{Q(\X)Q(\widetilde{\X}_{\sim u},\X_u)} = \Esmall{Q(\X_{\sim
    u},\X_u)Q(\widetilde{\X}_{\sim u},\X_u)}$ then
\begin{equation*}
  S_u = g(\theta)
      = \frac{\theta_3-\theta_2^2}{\theta_1-\theta_2^2}.
\end{equation*}

The construction of an estimator of $S_u$ boils down to the
construction of estimators of the three quantities $\theta_1$,
$\theta_2$ and $\theta_3$. To do this, nested Monte Carlo sampling is widely-used.
Let $\{\X^{(i)}; i=1,\dots,n\}$ and
$\{\widetilde{\X}^{(i)}; i=1,\dots,n\}$ be independent Monte Carlo
samples from the law of $\X$. For each $i$, denote by $\X_u^{(i)}$ the
subvector of $\X^{(i)}$ whose components are those of $\X^{(i)}$
indexed by $u$. Likewise, denote by $\X_{\sim u}^{(i)}$ the subvector
of $\X^{(i)}$ whose components are those of $\X^{(i)}$ not indexed by
$u$, and denote by $\widetilde{\X}_{\sim u}^{(i)}$ the subvector of
$\widetilde \X^{(i)}$ whose components are those of
$\widetilde \X^{(i)}$ not indexed by $u$.  An estimator of $S_u$ is given by
\begin{equation}\label{eq:sensindexestimator}
  \widehat S_u = g(\widehat\theta) =
  \frac{\widehat\theta_3-\widehat{\theta}_2^2}{\widehat\theta_1-\widehat{\theta}_2^2}
\end{equation}
where
\begin{equation}\label{eq:estimators}
  \left.
    \begin{array}{rcl}
      \widehat{\theta}_1 &=& \frac{1}{n}\sum_{i=1}^n\widehat{Q}_m(\X^{(i)})^2\\
      \widehat{\theta}_2 &=& \frac{1}{n}\sum_{i=1}^n\widehat{Q}_m(\X^{(i)})\\
      \widehat{\theta}_3 &=&  \frac{1}{n}\sum_{i=1}^n%
                             \widehat{Q}_m(\X^{(i)})%
                             \widetilde{Q}_m(\widetilde{\X}^{(i)}_{\sim u},\X^{(i)}_{u})\\
    \end{array}\right\}
\end{equation}
and
\begin{eqnarray*}
\widehat{Q}_m(\X^{(i)}) &=&
   \frac{1}{m}\sum_{k=1}^m\varphi\left(\X^{(i)},Z^{(i,k)}\right)\\
\widetilde{Q}_m(\widetilde{\X}^{(i)}_{\sim u},\X^{(i)}_{u}) &=&
   \frac{1}{m}\sum_{k=1}^m\varphi(\widetilde{\X}_{\sim u}^{(i)},\X_u^{(i)},%
                                                \widetilde Z^{(i,k)});
\end{eqnarray*}
here the objects $\{Z^{(i,k)},\widetilde{Z}^{(i,k)}; k=1,\cdots
,m; i=1,\cdots ,n\}$, are independent and identically distributed random
variables, independent of $\{\X^{(i)},\widetilde \X^{(i)};
i=1,\dots,n\}$, representing the randomness of the user's model. 
If $m$ is fixed and $n\to\infty$, then
\begin{equation*}
  \sqrt{n}\left(
    \widehat{S}_u-S_u
      \left[
        1 -
         \frac{
           \expec\var(\varphi(\X,Z)|\X)
         }{
           \expec\var(\varphi(\X,Z)|\X) 
              + m \var\expec(\varphi(\X,Z)|\X)
           }
      \right]
  \right) 
\end{equation*}
converges to a centered normal distribution with some variance
$\sigma_m^2$ depending on $m$. If, moreover, $m\to\infty$ such that
$\sqrt{n}/m\to 0$ then $\sqrt{n}(\widehat{S}_u-S_u)$ goes to a
centered normal distribution with variance
$\lim_{m\to\infty}\sigma_m^2$. For more details, see~\cite{mazo}. A
theoretically-guided choice for $m$ and $n$ that penalizes bad
rankings of the sensitivity index estimates
$\widehat S_1,\dots,\widehat S_p$ was given in~\cite{mazo}.


\section{Nested Monte Carlo estimation}
\label{sec:setting}
We say that an estimator
$ \widehat\theta$ is a nested Monte-Carlo estimator if
\begin{equation}
  \widehat\theta=\frac{1}{n}\sum_{i=1}^n
  \phi_1\left(\frac{1}{m}\sum_{k=1}^m\varphi_1(\xi^{(i)},\zeta^{(i,k)})\right) \cdots
  \phi_p\left(\frac{1}{m}\sum_{k=1}^m\varphi_p(\xi^{(i)},\zeta^{(i,k)})\right), \label{nmc}  
\end{equation}
where $\phi_1,\varphi_1,\dots,\phi_p,\varphi_p$ are measurable
functions and the sets of random vectors
$\{\xi^{(i)},\zeta^{(i,k)},\,k=1,\dots,m\}$, $i=1,\dots,n$, are
independent and identically distributed.  The estimator~\eqref{nmc} is
called a nested Monte Carlo estimator because it is built from two
nested Monte Carlo loops: an outer loop is used to simulate the random
vectors $\xi^{(1)},\dots,\xi^{(n)}$, and, for each $\xi^{(i)}$, an
inner loop is used to simulate $\zeta^{(i,k)}$,
$k=1,\dots,m$. The random vectors $\{\zeta^{(i,k)}\}$ and
$\{\xi^{(i)}\}$ are assumed to be independent. With this definition,   each of the
three estimators
$\widehat{\theta}_1,\widehat{\theta}_2,\widehat{\theta}_3$
in~\eqref{eq:estimators} is a nested Monte Carlo estimator: the outer
loop simulates the ``explorations'' $\X^{(i)},\widetilde{\X}^{(i)}$,
$i=1,\dots,n$, and the inner loop simulates the ``repetitions''
$\varphi(\X^{(i)},Z^{(i,k)}),\varphi(\widetilde{\X}_{\sim
  u}^{(i)},\X_u^{(i)},\widetilde Z^{(i,k)})$, $k=1,\dots,m$.

If, for each
$i=1,\dots,n$ and conditionally on $\xi^{(i)}$, the random variables
$\zeta^{(i,k)}$ are independent and identically distributed, and if
\begin{equation}
  \theta:=\mathbb{E}\left(
  \phi_1\left(\Ec{\varphi_1(\xi^{(1)},\zeta^{(1,1)})}{
      \xi^{(1)}}\right)\cdots
  \phi_p\left(\Ec{\varphi_p(\xi^{(1)},\zeta^{(1,1)})}{
      \xi^{(1)}}\right)\right), \label{eq:theta}
\end{equation}
then we have
\begin{align} \label{mse_theta}
  \mathbb{E}\,\| \widehat\theta-\theta\|^2
  &= \frac{1}{n}\text{Trace}\left(\Sigma_m\right)+ \|b_m\|^2,
\end{align} 
where $\Sigma_m$ is the variance-covariance matrix of
$$\phi_1\left(\frac{1}{m}\sum_{k=1}^m\varphi_1(\xi^{(1)},\zeta^{(1,k)})\right) \cdots
  \phi_p\left(\frac{1}{m}\sum_{k=1}^m\varphi_p(\xi^{(1)},\zeta^{(1,k)})\right).$$

The MSE in Equation \eqref{mse_theta}  is a sum of two terms: a
variance term which involves the ratio $1/n$ and a bias term that is
function of $m$  only.  

\begin{assumption}
 $b_m\to 0$ and $\Sigma_m\to \Sigma$ as $m\to +\infty$ \label{ass:conv_b}
\end{assumption}

 Assumption~\ref{ass:conv_b} is
fulfilled by Nested Monte Carlo estimators provided that the functions
$\phi_1,\dots,\phi_p$ in~\eqref{nmc} have good properties such as
boundedness, Lipschitz continuity or  boundedness of derivatives
\citep{giorgi,giorgi_weak_2020}. In the context of global sensitivity analysis for stochastic model,  \cite{mazo}
showed that Assumption \ref{ass:conv_b} is satisfied by $\widehat{\theta}_1,\widehat{\theta}_2,\widehat{\theta}_3$ (see Equation \eqref{eq:estimators})  associated with Sobol' index estimators.

If Assumption~\ref{ass:conv_b} is satisfied, it holds that
$\lim_{n,m\to +\infty}\mathbb{E}\,\| \widehat\theta-\theta\|^2=0$ and
thereby $ \widehat\theta$ converges in quadratic mean to $\theta$.
If, moreover, the rate at which $b_m$ goes to zero is known, then
optimal rates for $n=T^{1-\eta}$ and $m=T^{\eta}$, $\eta\in [0,1)$, can be calculated in terms of the total
number of runs $T$. 
For instance, in the case $p=1$ and $\phi_1$ is a smooth function with
uniformly bounded third derivative, \cite{juneja} showed that
$b_m=O(1/m)$ and thus $  \mathbb{E}\,\|
\widehat\theta-\theta\|^2=O(1/n+1/m^2)$, leading to $O(T^{2/3})$
reached for $n=O(T^{2/3})$ and $m=O(T^{1/3})$.  If $\phi_1$ is simply
Lipschitz continuous, \cite{rainforth} rather found that $
\mathbb{E}\,\| \widehat\theta-\theta\|^2=O(1/n+1/m)$, leading to a MSE
of order $O(T^{1/2})$, ensured by the choice of $n=O(T^{1/2})$ and $m=O(T^{1/2})$.

\subsection{Functions of nested Monte Carlo estimators}\label{sec:functionsNMCestimators}

Let $g: \mathcal{D}\rightarrow \mathbb{R}$ be a continuous,
non-constant function and let $\widehat{\theta}$ be as
in~(\ref{eq:estimators}). Note that the estimator $g( \widehat\theta)$
converges in probability to $g(\theta)$ but it does not necessarily
hold that its MSE $\mathbb{E}\,(g( \widehat\theta)-g(\theta))^2$
converges to 0.  In fact, the MSE may even be infinite, depending on
the properties of $g$ and the law of $\widehat\theta$. Since the law
of $\widehat{\theta}$ is unknown, it is important to find verifiable
conditions on $g$ under which the MSE is kept under control for all
possible laws of $\widehat{\theta}$.  For instance, if $g$ is Lipschitz continuous, then
there exists an constant $L$ such that
$$
\mathbb{E}\,\left(g( \widehat\theta)-g(\theta)\right)^2\leq L^2\mathbb{E}\,\| \widehat\theta-\theta\|^2.
$$

This bound could be minimized as in~\citet{juneja}, see above.
If $g$ is assumed to be a twice continuously differentiable such that
its Hessian matrix denoted $\nabla^2g$ is uniformly bounded, then, up
to existence of some moments of $ \widehat\theta$, the combination of
a Taylor-Lagrange expansion and convexity inequalities would yield
\begin{multline}\label{eq:IneqHessianUniformlyBounded}
    \mathbb{E}\,\left(g( \widehat\theta)-g(\theta)\right)^2\leq 4\mathbb{E}\left(\nabla g(\mu_m)^\top\left( \widehat\theta-\mu_m\right)\right)^2+ 2\left(g(\mu_m)-g(\theta)\right)^2\\
    + L'^2\left(\mathbb{E}\|\widehat{\theta} - \mu_m\|^4\right),
  \end{multline}
  with $L'>0$ such that $\sup_{x\in \mathcal{D}}\|\nabla^2g(x) \|_F\leq L'$ and $\|\cdot\|_F$ denotes the Frobenius norm.
From this upper bound, we could derive another bound by noticing that $\E{\|\widehat{\theta}-\mu_m\|^4}=O(1/n^2)$. This follows from  the Marcinkiewicz–Zygmund inequality \citep{zygmund}, provided that suitable moments exist.
It would then hold that:
\begin{align}
     \mathbb{E}\,\left(g( \widehat\theta)-g(\theta)\right)^2\leq 4\mathbb{E}\left(\nabla g(\mu_m)^\top\left( \widehat\theta-\mu_m\right)\right)^2+2\left(g(\mu_m)-g(\theta)\right)^2+O(1/n^2).\label{approx-decomp}
\end{align}
Equation \eqref{approx-decomp} highlights three terms in the
right-hand side: the first term
stands for a variance term, the second term represents a squared bias
term and the third term is a negligible approximation error. 

However, in practice,  in many interesting applications, the function $g$ does not have good enough
properties such as above.  For instance,  uniform boundedness of derivatives of $g$ is a very strong condition in general.  In order to weaken such a condition,  consider  the following one:
\begin{equation}
  \label{assumption_2}
\E{ \sup_{\lambda\in [0,1]}
\|\nabla^2 g\left(\lambda\widehat{\theta}+(1-\lambda)\mu_m\right)\|_F^4}=O(1) \quad \text{ as } n,m \to \infty.
\end{equation}

Under the condition~\eqref{assumption_2},  the  decomposition \eqref{approx-decomp}   holds but asymptotically when both $n$ and $m$ go to infinity.
However, whether or not this condition is true depends on the law of
$\widehat{\theta}$, and hence hard to check in practice. For
instance, 
in the case of Sobol' index estimators defined in Section
\ref{sec:context},  this comes down to require that
$\mathbb{E}[(\widehat{\theta}_1-\widehat{\theta}_2^2)^{-\alpha}] $
with $\alpha >0$ exists, whereas the probability distribution of
$\widehat{\theta}_1-\widehat{\theta}_2^2$ is unknown.

In the next section,  we propose  a condition, weaker
than the one in Equation \eqref{assumption_2}, that does not put
constraints  on the distribution of $\widehat\theta$.

 \subsection{A regularization method to address functions with unbounded Hessian matrix}\label{sec:regularization}

 The idea is to introduce a "slight perturbation" $g_h$ (with $h\in (0,1)$) of the function $g$ that approaches the true $g$ as $h\to 0$ and such that~\eqref{assumption_2} is satisfied with $g_h$ in place of $g$.
 Introducing  $g_h$ can be thought as a way to "transport" the original estimator $ \widehat\theta$ to regions of $\mathcal{D}$ where control of moments of $g_h( \widehat\theta)$ is possible without additional conditions on the law of $ \widehat\theta$. The goal is to ``get away''
  from certain regions of the parameter space where the Hessian of $g$
may explode. 
For this, consider a family of functions $\{g_h:\mathcal{D}\rightarrow \mathbb{R}, h\in (0,1)\}$ such that for all $h$, $g_h$ is twice continuously differentiable and for all $x\in \mathcal{D}$, $\lim_{h\to 0}g_h(x)=g(x)$. Moreover, assume $g_h$ satisfies the following assumption:

\begin{assumption}
There exists a constant $C$ independent of $h$ such that, for all $h \in (0,1)$:
$$
\limsup_{n,m\rightarrow \infty}\mathbb{E}\left(\sup_{\lambda\in [0,1]}\|\nabla^2 g_h(\lambda\widehat\theta+(1-\lambda)\mu_m)\|_F^4\right)\leq C.
$$
\label{assumption_3}
\end{assumption}

\begin{remark}
  Let us note that for some
  interesting classes of functions $g_h$ and estimators
  $\widehat{\theta}$, Assumption~\ref{assumption_3} can be checked
  without any knowledge of the law of $\widehat{\theta}$ but the
  existence of a sufficiently large number of finite moments. Indeed,
  in view of the linearity property of the expectation, this happens,
  for instance, if $\Vert g_h(\widehat{\theta}) \Vert_{\mathrm{F}}^2$
  is bounded by a polynomial in $\widehat{\theta}$ and
  $\widehat{\theta}$ is an empirical average of data with finite
  moments. That the Hessian be bounded in the Fr\"obenius norm is
  precisely what we would like to achieve by regularization. In other
  words, the latter property may not hold for the function $g$ while
  it holds for the function $g_h$.
\end{remark}

The advantage of having such a family of functions is that $\mathbb E(g_h(\widehat{\theta}) - g_h(\theta))^2$, the ``perturbed MSE'', could be bounded with an approximate upper bound in the form of Equation~\eqref{approx-decomp} with $g=g_h$.  A direct implication of this is that for $h$ fixed,  $ \lim_{n,m\to\infty}\mathbb E\left( g_h(\widehat\theta) - g_h(\theta) \right)^2=0$.
Further, if the regularized estimator $g_h(\widehat{\theta}) $ is used
for estimation of $g(\theta)$ instead of  $g(\widehat{\theta})$, the
corresponding mean squared error is bounded in Proposition~\ref{thm:basic-inequality}.

\begin{proposition}\label{thm:basic-inequality}
  Under Assumptions~\ref{ass:conv_b} and \ref{assumption_3}, for every $h\in (0,1)$:
 \begin{align} \label{inequality-gn}
    &\E{g_h(\widehat{\theta})-g(\theta)}^2 \le 4\left(1+p_{n,m}(h)\right)\left( V_{n,m}(h) + B_m(h)^2\right),
  \end{align}
  where $\limsup_{n,m\to \infty}p_{n,m}(h)=0$ and $V_{n,m}(h):=\mathbb{E}(\nabla
              g_h\left(\mu_m\right)^{\top}
              (\widehat{\theta}-\mu_m))^2$ and $B_m(h):=g_h\left(\mu_m\right)-g(\theta)$.
\end{proposition}

Remark that $B_m(h)=(g_h(\mu_m)-g_h(\theta))+(g_h(\theta)-g(\theta))$ accounts for two types of bias. The first term in the sum stands for the bias of the estimator $\widehat{\theta}$; it vanishes as soon as $\mu_m\to\theta$. The second term is a regularization error that stems from the use of $g_h$ instead of $g$. This term does not vanish asymptotically in general. However, it goes to zero as $h\to 0$, yielding Corollary~\ref{cor:hToZero}.
\begin{corollary}
  \label{cor:hToZero}
  We have $\lim_{h\to 0}\lim_{n,m\to\infty} \E{g_h(\widehat{\theta})-g(\theta)}^2 = 0$.
\end{corollary}

Note that the limits in Corollary~\ref{cor:hToZero} cannot be
interchanged. Additional conditions on  the gap $\mathbb{E}(g(\widehat{\theta})-g_h(\widehat{\theta}))^2$ would be needed to get the convergence of the MSE of $g(\widehat{\theta})$.


\section{Application to Sobol' index estimators}
\label{sec:application-to-sa}
Let 
$\widehat{\theta}=(\widehat{\theta}_1,\widehat{\theta}_2,\widehat{\theta}_3)$
be as in~\eqref{eq:estimators}.
Recall that
$\theta_1 = \Esmall{Q(\X)^2}$,
$\theta_2 = \Esmall{Q(\X)}$,
$\theta_3=\Esmall{\Ec{Q(\X)}{\X_u}^2} =
\Esmall{Q(\X)Q(\widetilde{\X}_{\sim u},\X_u)}$, 
so that 
$\mu_m = (\mu_{m1},\mu_{m2},\mu_{m3}) = (\theta_1+b_{m1},\theta_2+b_{m2},\theta_3+b_{m3})$ and 
$b_m = (b_{m1},b_{m2},b_{m3}) = 
(\mathbb E\Vc{\varphi(\X,Z)}{\X}/m,0,0)$.
Recall that $\widehat{S}_u = g(\widehat{\theta})$, where the function \begin{equation}
  g : (x_1,x_2,x_3) \mapsto
  (x_3-x_2^2)/(x_1-x_2^2)\label{eq:THE_g}\end{equation} is a
twice-continuously differentiable function over its definition domain.
But unfortunately, $g$ is unbounded. The form of such a function makes
the study of the MSE of  $\widehat{S}_u$ difficult unless strong
conditions are imposed on the output distribution of the stochastic
model. This could explain why, in spite of our efforts, it is difficult
to find in the literature material on the convergence in quadratic
mean of such estimators. The approach introduced in Section
\ref{sec:setting} allows to bypass the unboundedness issue. We shall
apply the results of Section~\ref{sec:setting} to bound the MSE of the
regularized version of the Sobol' index estimator and propose an allocation method. 
  
Throughout this section,  it is assumed that $\E{Q(\X)^{16}}<+\infty$.  We consider the family of function $g_h(x):=g(x+h\mathbf{u})$ with $\mathbf u=(1,0,0)$.
In order to provide an upper bound for $\mathbb{E}(g_h(\widehat{\theta})-g(\theta))^2$ as in Proposition~\ref{thm:basic-inequality},
it is necessary to fulfill Assumption \ref{ass:conv_b} and \ref{assumption_3}. Assumption~\ref{ass:conv_b} is obviously satisfied. 

\begin{theorem}\label{res:main}
 With $\widehat{\theta}$, $g_h$ and $g$ as above,
 Assumption~\ref{assumption_3} holds.  There is a constant $\bar{C}$
 independent of $h$ such that, for every $h \in (0,1)$,
    \begin{align*}
    \E{g_h(\widehat{\theta})-g(\theta)}^2
      \leq
     \bar{C}\left(
       \frac{1}{n}+\left(\frac{\E{\Vc{\varphi(\X,Z)}{\X}}}{m}
     \right)^2 + h^2\right),
  \end{align*}
  as $n,m\to \infty$.
\end{theorem}

\begin{corollary}\label{cor:convrate}
  With $\widehat{\theta}$, $g_h$ and $g$ as above, we have
  $\lim_{h\to 0}\lim_{n,m\to\infty}\mathbb
  E(g_h(\widehat{\theta})-g(\theta))^2 = 0$.
\end{corollary}
As $n,m\to +\infty$,  it occurs that the smaller $h$ is, the smaller the MSE gets.  Therefore, in practice, parameter $h$ should be taken in order to approximate the true $g(\theta)$.

%

To balance the trade-off between the number of repetitions and
explorations, we proceed with a heuristic. For $h$ small enough,
one could assume that $h^2$ is
negligible with respect to the quantity
$\text{BVT}:=1/n+(\E{\Vc{\varphi(\X,Z)}{\X}}/m)^2$.
This is, therefore, the quantity we shall be interested in for
optimizing the numbers of repetitions and explorations. Let $T=mn$ be
the number of available samples to compute the Sobol' index
estimator. From nested Monte Carlo theory (see beginning of
Section~\ref{sec:setting}), we already know that 
$m$ and $n$ must be of order $T^{1/3}$ and $T^{2/3}$ respectively so as to obtain optimal convergence rate.


However, an asymptotic order is not a specific value. We expect that the greater the variance of intrinsic randomness, the larger the value of optimal $m$ should be. Thus, in this section, the goal consists in proposing a value of $m$ that adapts to the intrinsic randomness of the stochastic model. 

Under the constraint $nm=T$, the optimal convergence rate of the BVT is obtained when $n_{opt}$ is of order $T^{2/3}$ and $m_{opt}$ is of order $T^{1/3}$. Let $m_{opt}=\kappa T^{1/3}$ where $\kappa>0$. Then,  $n_{opt}=\kappa^{-1}T^{2/3}$.  Thus:
$$
\mathrm{BVT}=\kappa T^{-2/3}+\frac{(\mathbb{E}\Vc{\varphi(\mathbf{X},Z)}{\X})^2}{\kappa^2}T^{-2/3}.
$$
Coefficient $\kappa$ can be chosen such that $\kappa T^{-2/3}+\frac{(\mathbb{E}\Vc{\varphi(\mathbf{X},Z)}{\X})^2}{\kappa^2}T^{-2/3}$ is the smallest over $\kappa>0$. The minimum of such a quantity is reached at  $\kappa_{opt}=\left(2~ (\mathbb{E}\Vc{\varphi(\mathbf{X},Z)}{\X})^2\right)^{1/3}$.  Therefore:

\begin{equation}\label{eq:mopt}
  m_{opt}(T)=\left(2~(\mathbb{E}\Vc{\varphi(\X,Z)}{\X})^2\right)^{1/3}T^{1/3}.
\end{equation}
Therefore, the number of repetitions suggested above ensures first
that the BVT converges at optimal rate and second adapts to the model's intrinsic randomness because that $m_{opt}(T)$ depends on 
$\mathbb{E}\Vc{\varphi(\X,Z)}{\X}$, which 
quantifies the part of the total variance $\V{\varphi(\X,Z)}$ that is not attributed to the inputs $\X$; and so, that measures the influence of the intrinsic noise of the stochastic model $\varphi(\X,Z)$.
Finally, it also appears that $m_{opt}(T)$ depends on both $T$ and the function $\varphi$.  The dependence with respect to $T$ guarantees that $m_{opt}(T)$ grows as $T$ gets large. Besides, the dependence with respect to $\varphi$ means that even if $m_{opt}(T)$ remains proportional to $T^{1/3}$, it also varies with respect to the chosen QoI of the stochastic model $f$.


\section{An algorithm of estimation of Sobol' indices for stochastic models}
\label{algo-illus-section}
\subsection{Algorithms}

This section is devoted to the practical implementation of the
bias-variance trade-off strategy when performing SA for some QoI of a stochastic model. Recall that $f$ is a stochastic model as in~(\ref{stoch-model}) and we are interested in carrying out SA of a QoI under the form~(\ref{qoi-form}), that is, $Q(\X)=\Ec{\varphi\left(\X,Z\right)}{\X}$ in order to measure the impact of $l$ groups of inputs $u_i\subset\{1,\dots,p\}$, $i=1,\dots,l$. In other words, we are interested in estimating $S_{u_1},\dots,S_{u_l}$. We shall use at most $T\times (l+1)$ evaluations of $\varphi(\X,Z)$. Under the constraint $nm=T$, the  number of repetitions $m_{opt}$ found in~\eqref{eq:mopt} depends on  $\rho:=\mathbb{E}\Vc{\varphi(\X,Z)}{\X}$.  However, in practice, $\rho$ is often unknown. So, before sensitivity index estimation,  $\rho$ needs to be estimated. Notice that $\rho$ can be rewritten $\rho=\mathbb{E}(\varphi(\X,Z)-\varphi(\X,\widetilde{Z}))^2$, which suggests the following estimation procedure.

Consider $r_0$ i.i.d. samples of $\X$, denoted by $\X^{(1)}, \cdots , \X^{(r_0)}$,  and generate two outputs at each sample $\X^{(i)}$: $\left(\varphi(\X^{(1)},Z^{(1,1)}),\varphi(\X^{(1)},Z^{(1,2)})\right),
\cdots, $\\
$\left(\varphi(\X^{(r_0)},Z^{(r_0,1)}),\varphi(\X^{(r_0)},Z^{(r_0,2)})\right)$.
Thus:
$$
\widehat{\rho}=\frac{1}{r_0}\sum_{i=1}^{r_0}\left(\varphi(\X^{(i)},Z^{(i,1)})-\varphi(\X^{(i)},Z^{(i,2)})\right)^2
$$ 
is a consistent and unbiased estimator of $\rho$. It appears that the estimation of $\rho$ requires $2r_0$ evaluations of the model $\varphi(\X,Z)$. However, the maximal number of evaluations is $T\times (l+1)$. So, for index estimation procedure,  at most $T\times (l+1)-2r_0$ model evaluations are allowed. 

Therefore, our strategy consists in leveraging the model outputs used to estimate $\rho$ and then plugging and completing those outputs in order to compute sensitivity index estimates. This strategy relies on two algorithms: Algorithm \ref{algo-1} and Algorithm \ref{algo}.  Algorithm \ref{algo-1} enables to generate complementary outputs in addition to outputs already available after estimation of $\rho$. This allows to satisfy the constraint of the maximal number of model evaluations $T\times (l+1)$. This part helps to optimize the whole estimation procedure by using the model outputs already generated. Regarding Algorithm \ref{algo}, it effectively estimates indices in three steps based on pick-freeze procedure. First, it estimates $\rho$ and thereby compute $m_{opt}$ and $n_{opt}=T/m_{opt}$. Then, in the second step, by relying on Algorithm \ref{algo-1}, complementary outputs required for estimations are generated. In the final step, sensitivity index estimates are computed with respect to inputs or groups of inputs specified by the user.
\begin{algorithm}[H]
\caption{Completing model evaluations}
\hspace*{\algorithmicindent} \textbf{Inputs:} $n,m,\varphi, l,\left(\X^{(1)}, \cdots ,\X^{(T)}\right)$ \\
\hspace*{\algorithmicindent} \textbf{Data:} $\left(\varphi\left(\X^{(1)},Z^{(1,1)}\right),\varphi\left(\X^{(1)},Z^{(1,2)}\right)\right), \cdots ,\left(\varphi\left(\X^{(r_0)},Z^{(r_0,1)}\right),\varphi\left(\X^{(r_0)},Z^{(r_0,2)}\right)\right)$
\begin{algorithmic}
 \IF{$n\geq r_0$}
    \IF{$m \geq 2$}
   \FOR{$i=1,\cdots ,r_0$}
    \FOR{$k=3,\cdots ,m$}
   \STATE Compute $\varphi\left(\X^{(i)},Z^{(i,k)}\right)$ 
   \ENDFOR
   \ENDFOR
   \FOR{$i=r_0+1,\cdots ,n$}
    \FOR{$k=1,\cdots ,m$}
   \STATE Compute $\varphi\left(\X^{(i)},Z^{(i,k)}\right)$ 
   \ENDFOR
   \ENDFOR
    \ENDIF
     \IF{$m=1$}
   \FOR{$i=r_0+1,\cdots ,n-r_0-\lceil r_0/(l+1)\rceil$}
   \STATE Compute $\varphi\left(\X^{(i)},Z^{(i,1)}\right)$ 
   \ENDFOR
    \ENDIF
   \ENDIF
    \IF{$n < r_0$}
    \IF{$m > 2+ 2\lceil 1/(l+1)*(-1+r_0/n)\rceil$}
    \FOR{$i=1,\cdots ,n$}
    \FOR{$k=3,\cdots ,m- 2\lceil 1/(l+1)*(-1+r_0/n)\rceil$}
   \STATE Compute $\varphi\left(\X^{(i)},Z^{(i,k)}\right)$
                   
   \ENDFOR
   \ENDFOR
    \ENDIF
     \IF{$m \leq 2+ 2\lceil 1/(l+1)*(-1+r_0/n)\rceil$}
      \STATE  \textbf{Exit:} Budget already consumed
    \ENDIF
   \ENDIF
\end{algorithmic}
\label{algo-1}
\end{algorithm}

\begin{algorithm}[H]
  \caption{Estimation of Sobol' indices}
   \hspace*{\algorithmicindent} \textbf{Inputs:} $h, r_0, T, \varphi, w=\{u_1,\cdots ,u_l\}, \left(\X^{(1)}, \cdots ,\X^{(T)}\right),\left(\widetilde\X^{(1)}, \cdots ,\widetilde\X^{(T)}\right)$ 
  \begin{algorithmic}[1]
   \FOR{$i=1,\cdots ,r_0$}
    \FOR{$k=1,2$}
   \STATE Compute $\varphi\left(\X^{(i)},Z^{(i,k)}\right)$ 
   \ENDFOR
   \ENDFOR
   \STATE Compute $\widehat{\rho} \gets \frac{1}{r_0}\sum_{s=1}^{r_0}\left(\varphi\left(\X^{(i)},Z^{(i,1)}\right)-\varphi\left(\X^{(i)},Z^{(i,2)}\right)\right)^2$
   \STATE Compute $\widehat{m}_{opt}$ according to Equation \eqref{eq:mopt}
   \STATE $\widehat{n}_{opt} \gets \big\lfloor T/\widehat{m}_{opt}\big\rfloor$
   \STATE Run Algorithm \ref{algo-1}  with $m=\widehat{m}_{opt},n=\widehat{n}_{opt} $  to complete samples $\varphi\left(\X^{(1)},Z^{(1,1)}\right), \cdots ,\varphi\left(\X^{(r_0)},Z^{(r_0,2)}\right)$.
     \FOR{$j=1,\cdots ,l$}
   \FOR{$i=1,\cdots ,\widehat{n}_{opt}$}
    \FOR{$k=1,\cdots ,\widehat{m}_{opt}$}
   \STATE Compute  $\varphi\left((\widetilde{\X}_{\sim u_j}^{(i)},\X_{u_j}^{(i)}),\widetilde{Z}^{(i,k)}\right)$  
   \ENDFOR
   \ENDFOR
     \ENDFOR
   \STATE Compute sensitivity index estimates or their regularized estimates
 
    
  \end{algorithmic}
  \label{algo}
\end{algorithm}

Algorithm \ref{algo} requires: $h, r_0,T,\varphi, w$ and input samples.  The transformation $\varphi$ of the stochastic model is supplied as well as $w$ the set of inputs or groups of inputs whose indices are estimated.  In practice,  $r_0$, $T$ and $h$ must be chosen. We recommend to take $r_0$ small with respect to $T$ so as not to waste too much of the budget in the first stage of Algorithm \ref{algo}. Indeed, the estimator $\widehat{\rho}$ is unbiased and consistent and often provides good estimates even for small $r_0$.  Regarding $T$,  it  follows $T$ should be taken as large as possible depending on the computational cost of a run of  both $\varphi$ and the original model $f$. Furthermore, to ensure that the MSE has a precision $\varepsilon\in (0,1)$ with $h^2 \ll \varepsilon$,  $T$ must be roughly chosen larger than $\varepsilon^{-3/2}$ since the MSE is $O\left(T^{-2/3}\right)$. This provides approximations for practical choice of $T$. 

%
%

\section{Illustrations}
\label{illustration}
This section presents the performance of the estimators  of first-order and total indices computed by Algorithm \ref{algo} in the case of two toy stochastic models for which analytical values of indices are known: a linear model $f(X_1,X_2,Z)=1+ X_1+2X_2+\sigma Z$ with $\sigma >0$ and a stochastic version of the Ishigami function  $f^\prime(X_1,X_2,X_3,Z)=\sin X_1 +a\sin^2X_2+bX_3^4(\sin X_1)Z^2$ with $a, b>0$ (\cite{ishigami}).
For each value of $T=nm$, the estimators of Algorithm \ref{algo} are compared with two other arbitrary choices, namely, $(n,m)=(T/5,5)$ and $(n,m)=(T^{1/2},T^{1/2})$.


These two  choices above represent two different situations. The choice $(n,m)=(T/5,5)$ presents a case where the number of repetitions is constant and independent of $T$. This illustrates the situation where the bias does not get reduced so that it disturbs estimations no matter how large $T$ is. Regarding $(n,m)=(T^{1/2},T^{1/2})$, it shows that the case where the variance is not sufficiently reduced.

For illustrations, the product $T=mn$ is chosen in the set $T\in \{10^3,10^4,\cdots ,10^7\}$.  For each choice of the couple $(n,m)$, $N=100$ replications of estimations are carried out. The tuning parameter $r_0$ is set to $10$. 
Thus, $2r_0=20$ model evaluations are used to get the estimates $\widehat \rho$, $\widehat n_{opt}$, $\widehat m_{opt}$ in the first part of Algorithm~\ref{algo}, and then $T-2r_0\in\{10^3-20,10^4-20,\cdots ,10^7-20\}$ model evaluations are used to get the sensitivity index estimators with $n=\widehat n_{opt}$ and $m=\widehat m_{opt}$.
For both toy stochastic models, the QoI considered is the conditional expectation so that $\varphi=f$ or $\varphi=f^\prime$ depending on the model. For each value of $T$, the boxplots of sensitivity index estimates are plotted for each of the three choices.

\subsubsection*{Linear model}
Let  $f(X_1,X_2,Z)= 1+X_1+2X_2+\sigma Z$ where $\sigma >0$ and $X_1,X_2$ and $Z$ are i.i.d.  under the standard normal distribution.  Such model includes two uncertain parameters $X_1$ and $X_2$ with respective first-order Sobol' indices $S_1=1/5$ and $S_2=4/5$. 
Two values of $\sigma$ are considered:  $\sigma=1$  and $\sigma=5$.
 
\begin{figure}[H]
     \centering
     \begin{subfigure}[b]{\textwidth}
         \centering
         \includegraphics[scale=0.5]{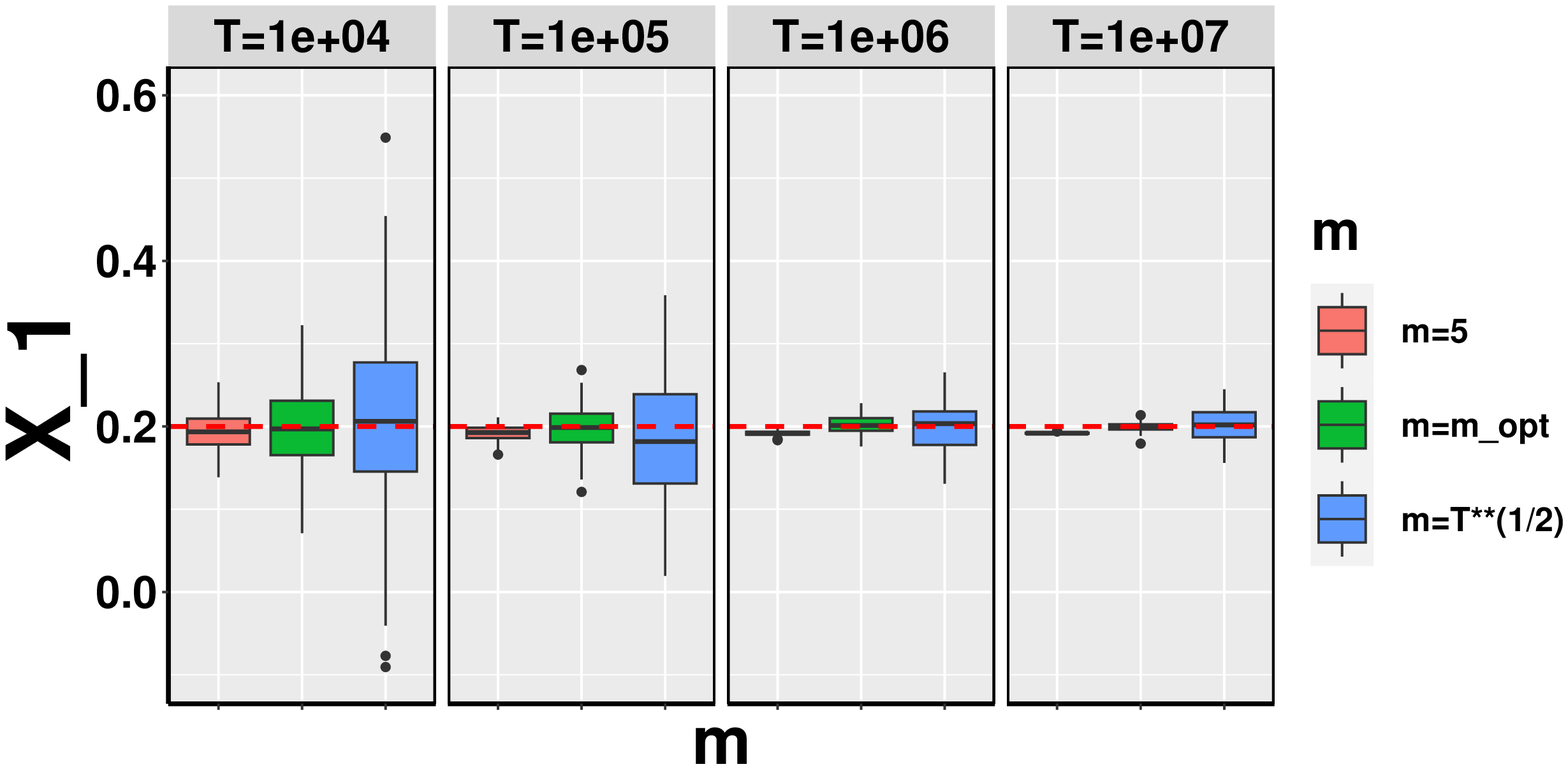}
         \label{fig:sig1x1}
     \end{subfigure}
    \\
     \begin{subfigure}[b]{\textwidth}
         \centering
         \includegraphics[scale=0.5]{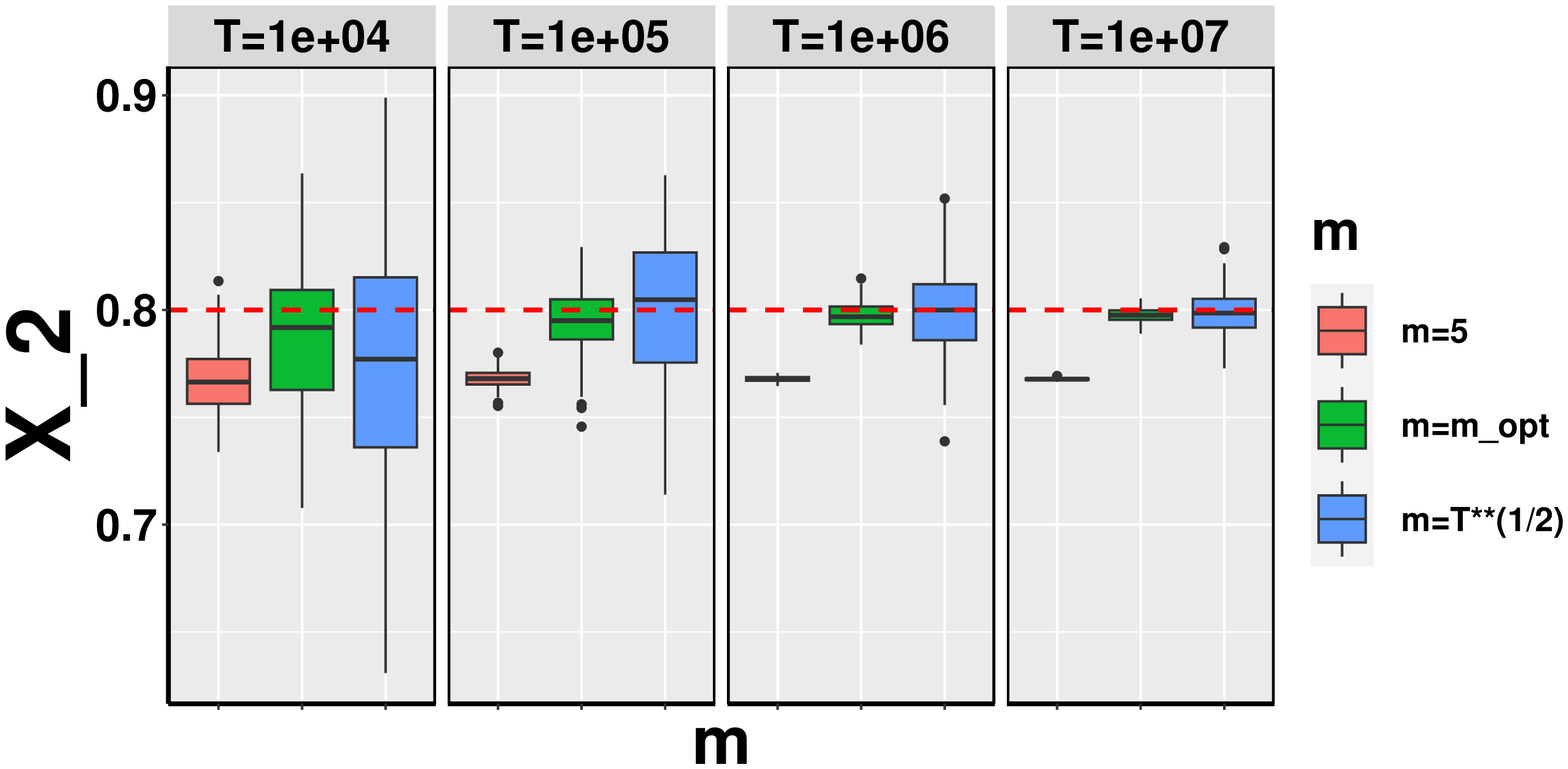}
         \label{fig:sig1x2}
     \end{subfigure}
        \caption{Boxplots of sensitivity index estimates of $X_1$ and $X_2$ with respect to  the linear model and for different values of $T$.  The standard deviation is $\sigma=1$ and regularization parameter $h=10^{-2}$. Three strategies of choice of $m$ are compared: $m=5$ (in red), $m=m_{opt}$ given by the trade-off strategy of Algorithm \ref{algo} (in green) and $m=T^{1/2}$ (in blue). The red dashed lines showed the true sensitivity index values, in particular, in this setting $S_1=0.2$ and $S_2=0.8$}
        \label{linear-mse}
\end{figure}

\begin{figure}[H]
     \centering
     \begin{subfigure}[b]{\textwidth}
         \centering
         \includegraphics[scale=0.5]{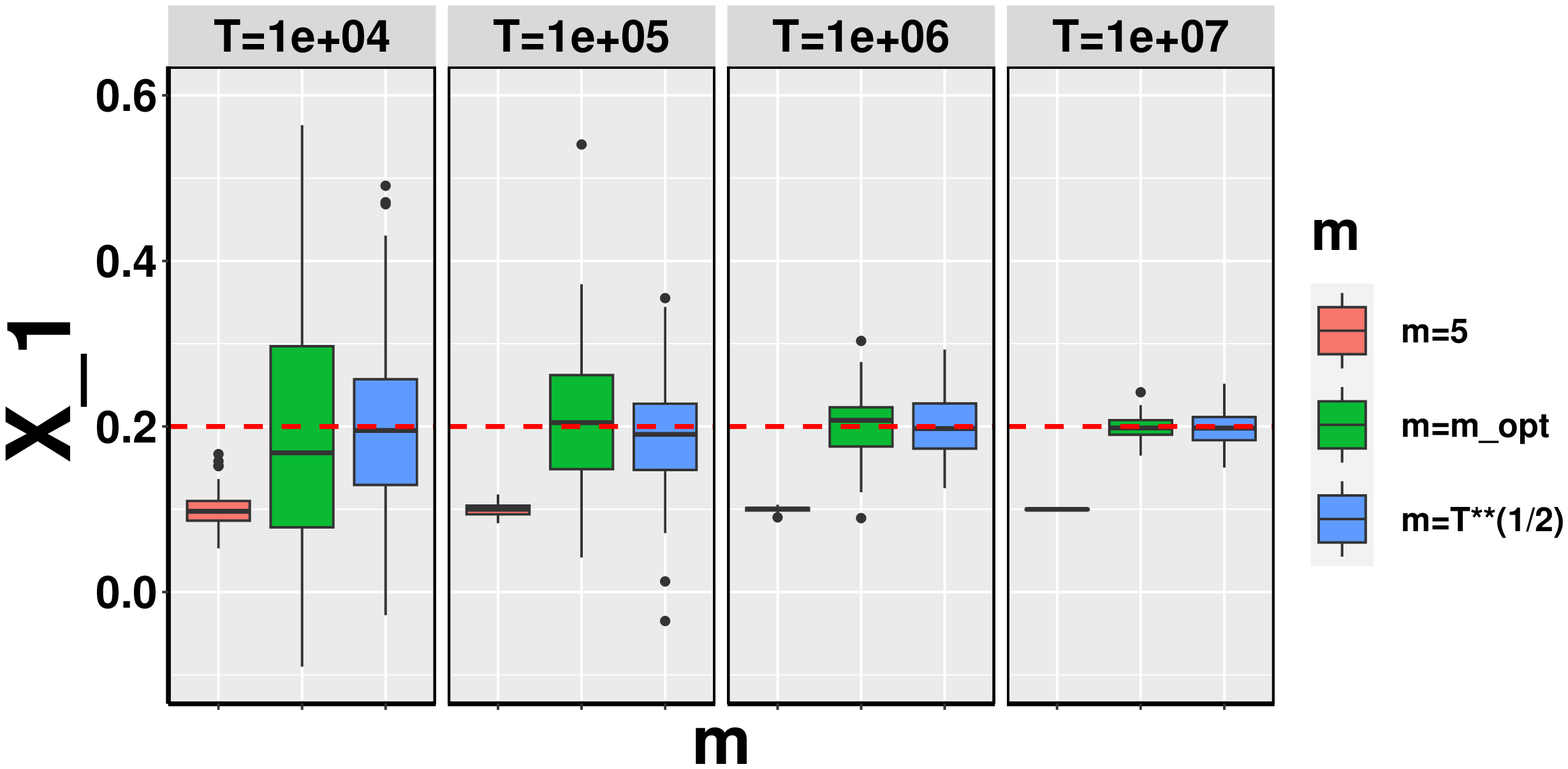}
         \label{fig:sig5x1}
     \end{subfigure}
    \\
     \begin{subfigure}[b]{\textwidth}
         \centering
         \includegraphics[scale=0.5]{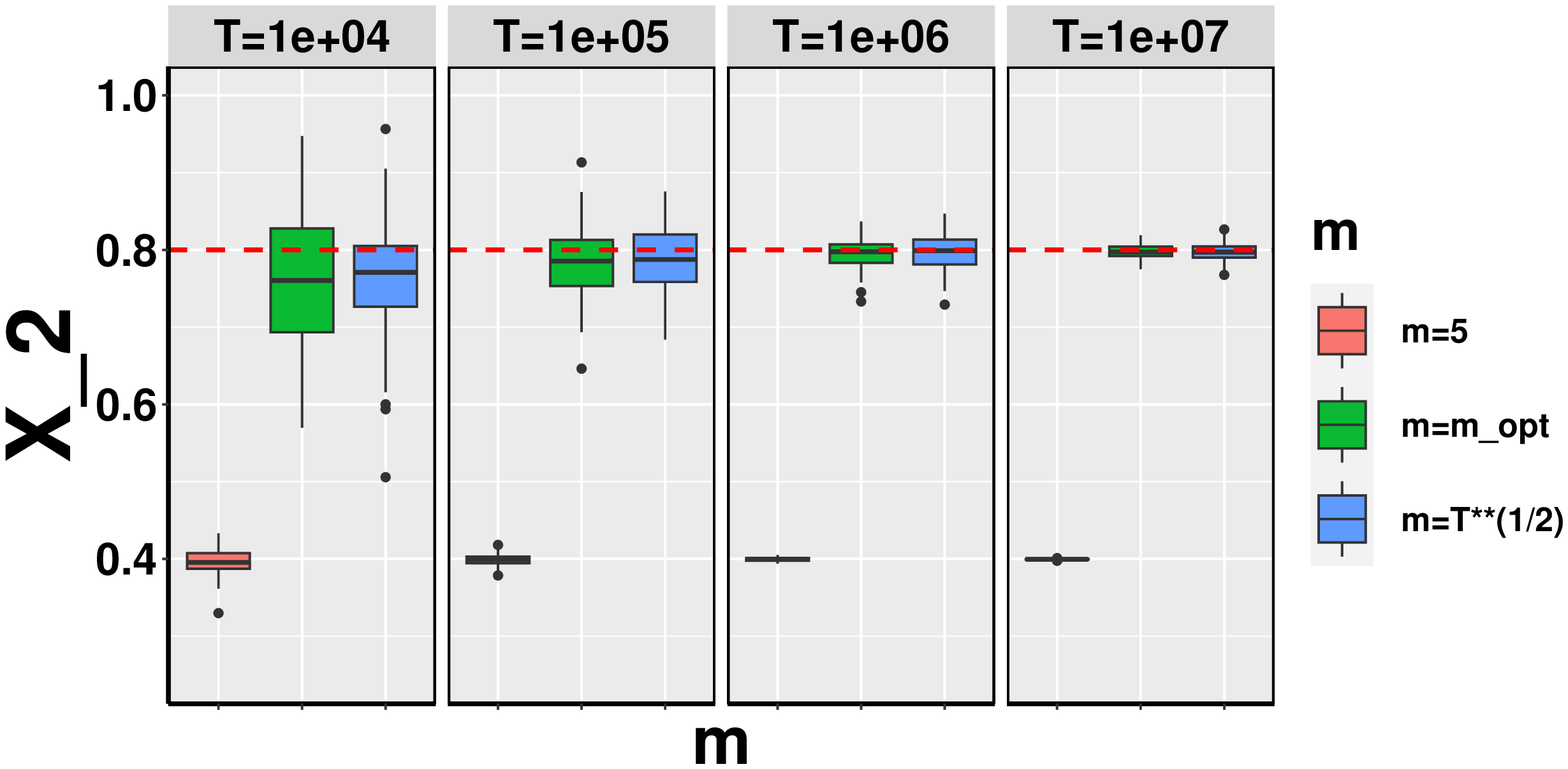}
         \label{fig:sig5x2}
     \end{subfigure}
        \caption{Boxplots of sensitivity index estimates of $X_1$ and $X_2$ with respect to  the linear model and for different values of $T$.  The standard deviation and regularization parameter are  $\sigma=5$ and $h=10^{-2}$.   Three strategies of choice of $m$ are compared: $m=5$ (in red), $m=m_{opt}$ given by the trade-off strategy of Algorithm \ref{algo} (in green) and $m=T^{1/2}$ (in blue). The red dashed lines showed the true sensitivity index values, in particular, in this setting $S_1=0.2$ and $S_2=0.8$}
        \label{linear-mse2}
\end{figure}

 Figure \ref{linear-mse} shows that the estimations obtained with Algorithm \ref{algo} are more efficient as $T$ gets large because both bias and variance are efficiently reduced.  Boxplots  highlight that the strategies $m=5$  and $m=T^{1/2}$ suffer respectively from bias and variance. Notice that in the case of the linear model under study,  $\mathbb{E}(\var(f(X_1,X_2,Z)\mid (X_1,X_2)))=2\sigma^2$; so the bias depends on $\sigma$. This explains why in the case $\sigma=5$ (see Figure \ref{linear-mse2}), even for large value of $T$, estimations resulting of the choice $m=5$ do not converge to their true targets.

\subsubsection*{A stochastic Ishigami function}

Let $f^\prime(X_1,X_2,X_3,Z)=\sin X_1 +a\sin^2X_2+bX_3^4(\sin X_1)Z^2$  such that with $a, b>0$, $X_1,X_2,X_3$ and $Z$ are independent with $X_1,X_2,X_3$ distributed under $\mathcal{U}\left([-\pi,\pi]\right)$ and $Z\sim \mathcal{N}\left(0,1\right)$.  The model $f^\prime$ is a modified version of benchmark function known as the Ishigami function in SA.  For this model, first-order Sobol' indices of inputs $X_1, X_2,$ and $X_3$ for the QoI $\mathbb{E}(f^\prime(X_1,X_2,X_3,Z)\mid X_1,X_2,X3)$ are respectively given by $S_1=\frac{1}{2}\frac{\left(1+\frac{b\pi^4}{5}\right)^2}{\frac{a^2}{8}+\frac{b\pi^4}{5}+\frac{b^2\pi^8}{18}+\frac{1}{2}}$, $S_2=\frac{\frac{a^2}{8}}{\frac{a^2}{8}+\frac{b\pi^4}{5}+\frac{b^2\pi^8}{18}+\frac{1}{2}}$ and $S_3=0$.  Parameters $a$ and $b$ are chosen with respect to \cite{levitan}: $a=7$, $b=0.05$ and   \cite{marrel2}: $a=7$, $b=0.1$.
\begin{figure}[H]
     \centering
     \begin{subfigure}[b]{\textwidth}
         \centering
         \includegraphics[scale=0.4]{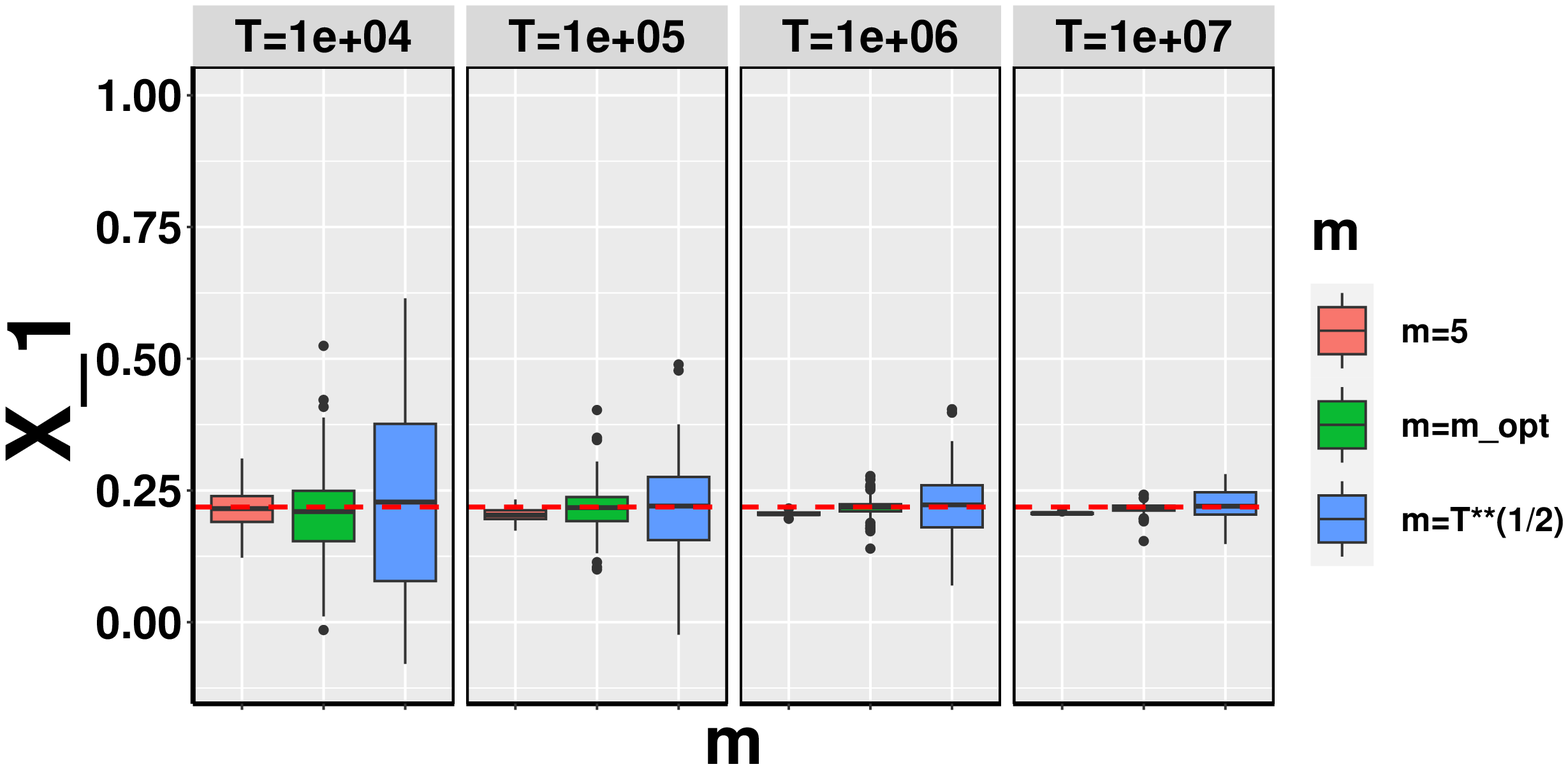}
         \label{fig:b1x1}
     \end{subfigure}
      \\
     \begin{subfigure}[b]{\textwidth}
         \centering
         \includegraphics[scale=0.4]{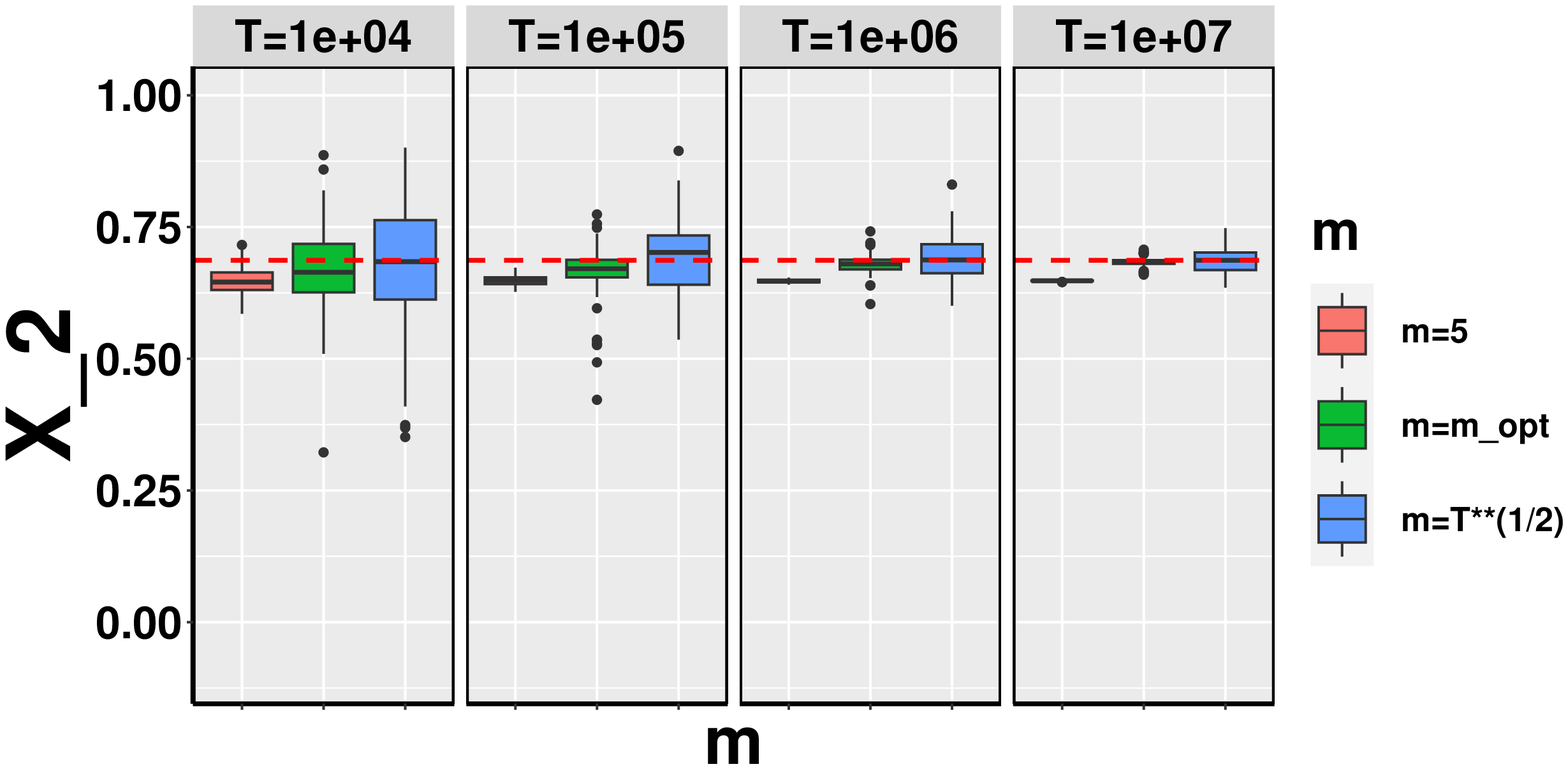}
         \label{fig:b1x2}
     \end{subfigure}
       \\
     \begin{subfigure}[b]{\textwidth}
         \centering
         \includegraphics[scale=0.4]{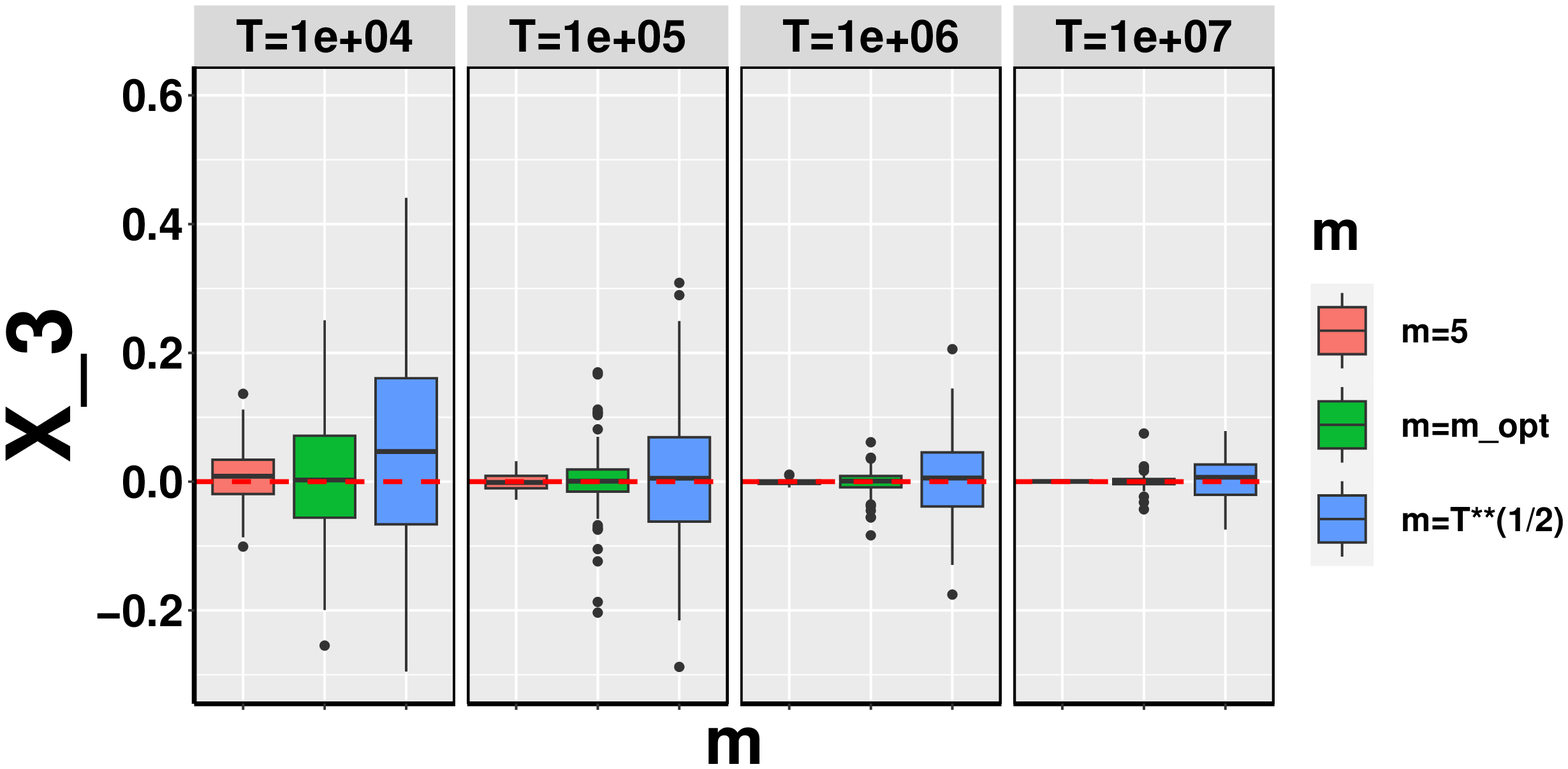}
         \label{fig:b1x3}
     \end{subfigure}
        \caption{Boxplots of sensitivity index estimates of $X_1$, $X_2$ and $X_3$ with respect to the stochastic version of Ishigami function (with $b=0.05$) for different values of $T$. Three strategies of choice of $m$ are compared: $m=5$ (in red), $m=m_{opt}$ given by the trade-off strategy of Algorithm \ref{algo} (in green) and $m=T^{1/2}$ (in blue). The red dashed lines showed the true sensitivity index values.}
        \label{ishigami-mse-005}
\end{figure}

\begin{figure}[H]
     \centering
     \begin{subfigure}[b]{\textwidth}
         \centering
         \includegraphics[scale=0.47]{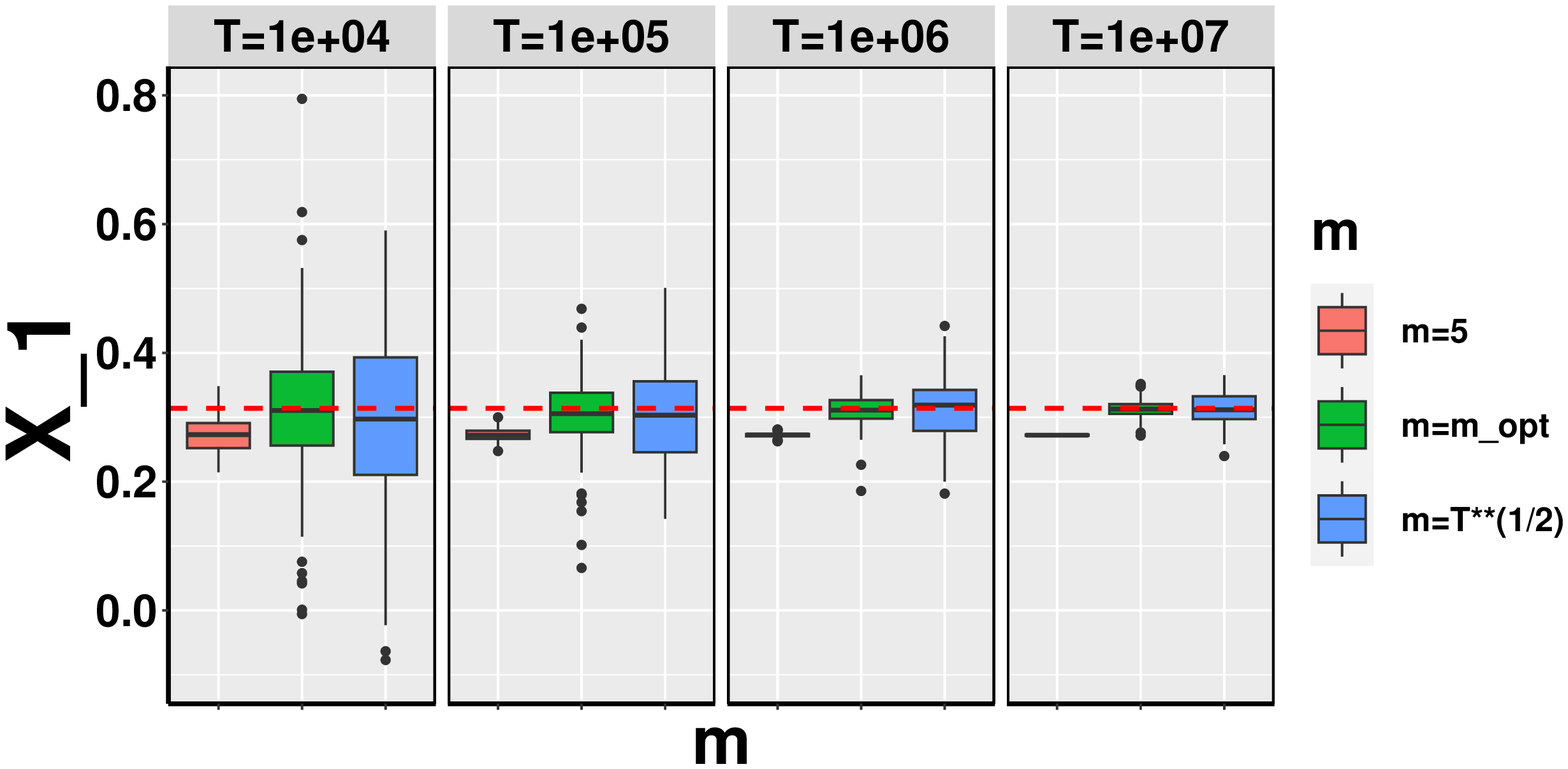}
         \label{fig:b2x1}
     \end{subfigure}
      \\
     \begin{subfigure}[b]{\textwidth}
         \centering
         \includegraphics[scale=0.47]{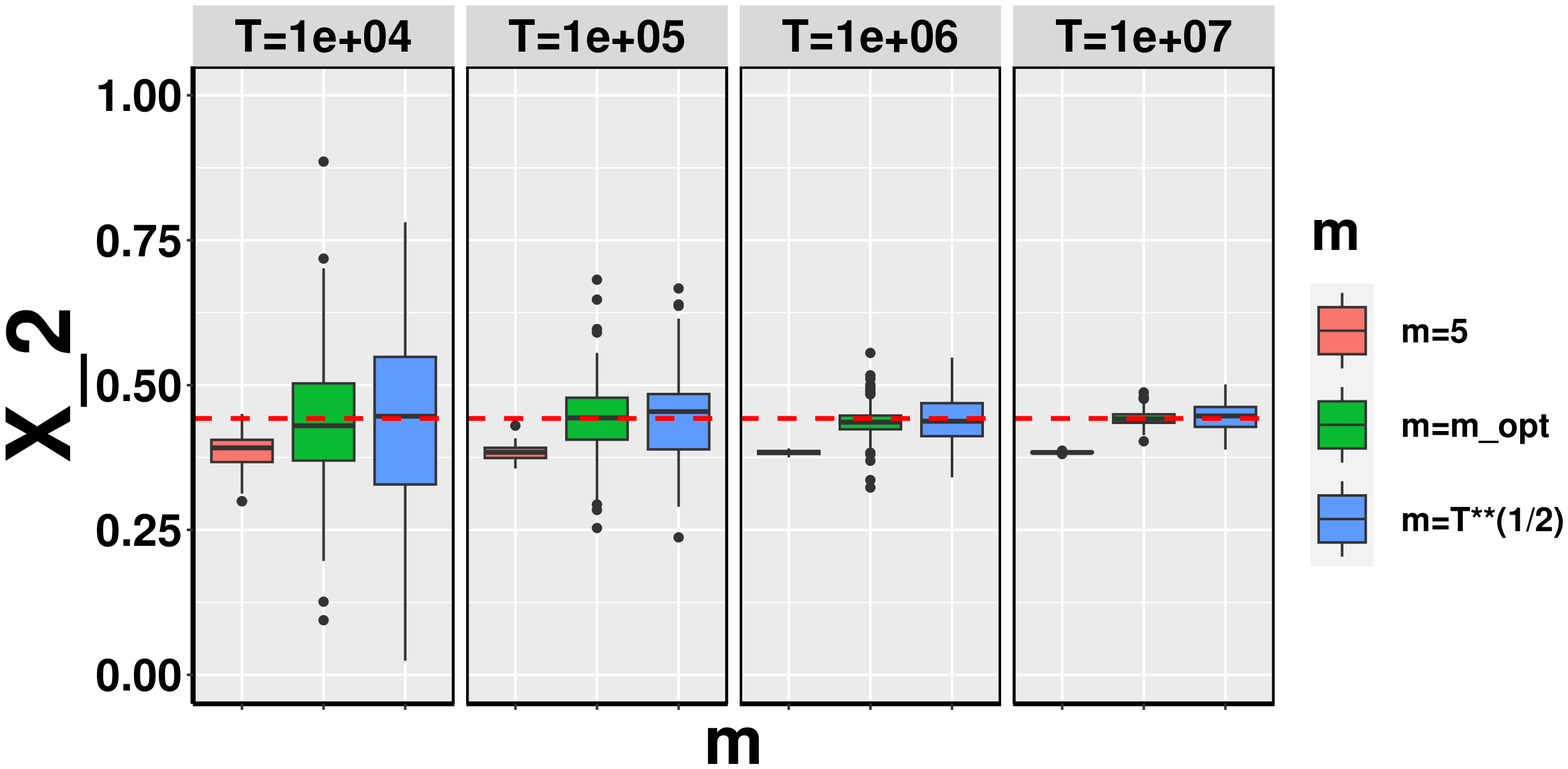}
         \label{fig:b2x2}
     \end{subfigure}
       \\
     \begin{subfigure}[b]{\textwidth}
         \centering
         \includegraphics[scale=0.47]{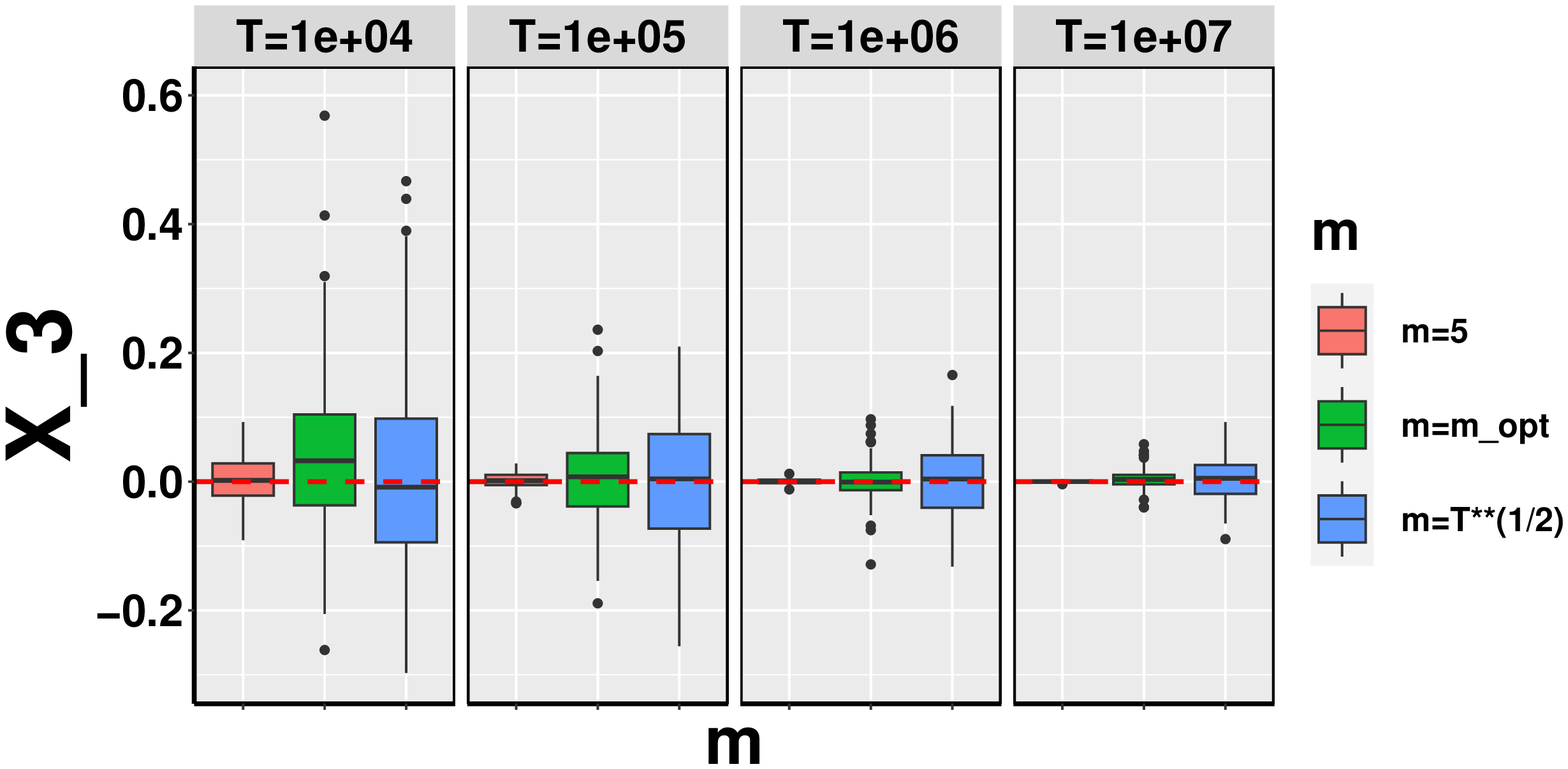}
         \label{fig:b2x3}
     \end{subfigure}
        \caption{Boxplots of sensitivity index estimates of $X_1$, $X_2$ and $X_3$ with respect to the stochastic version of Ishigami function (with $b=0.1$) for different values of $T$. Three strategies of choice of $m$ are compared: $m=5$ (in red), $m=m_{opt}$ given by the trade-off strategy of Algorithm \ref{algo} (in green) and $m=T^{1/2}$ (in blue). The red dashed lines showed the true sensitivity index values.}
        \label{ishigami-mse-01}
\end{figure}

Figures \ref{ishigami-mse-005} and  \ref{ishigami-mse-01} also reveal that estimations obtained by using Algorithm \ref{algo} are more efficient for large $T$.  Besides, remark that the term $bX_3^4\sin X_1$ that multiplies the intrinsic noise term $Z^2$ includes $b$ so that $\vert bX_3^4\sin X_1\vert \leq b\pi^4$. Then, $b$ allows to control the variance of the intrinsic noise term of the model.
In both cases, the strategy implemented in Algorithm \ref{algo} has better results.

Overall,  it appears that the strategy of Algorithm \ref{algo} provides better estimations  and estimators it generates converge faster.  In the particular case of $m=5$, it is noticeable that errors do not decrease when $T$ gets larger but rather they are quite constant. This is explained by the fact that the bias is constant since $m$ is constant. This illustrates the importance of varying the number of repetition when the total computational budget grows. Regarding the case $m=T^{1/2}$, it turns out that estimators do not converge with optimal rate compared to the case $m=m_{opt}$ due the variance part of the mean-squared error. Indeed, with $m=T^{1/2}$, the variance converges to $0$ at rate $T^{1/2}$ while the squared-bias converges at rate $T$. Then, the global convergence rate of the MSEs is $T^{1/2}$ that is slower than the rate $T^{2/3}$ of estimators built by Algorithm \ref{algo}.  These two cases clearly illustrate the bias-variance trade-off problem in Sobol' index estimation for stochastic models and they allow to show that the strategy proposed in this paper performs well.


\section{Conclusion}
\label{conclusion}

To balance the tradeoff between repetitions and explorations in the
estimation of Sobol' indices for stochastic models, a regularized
estimator was devised and its mean-squared error bounded, approximated
and minimized. A connection to functions of nested Monte Carlo
estimators was made, for which generic results were obtained,
especially useful when the Hessian matrix is unbounded. To estimate
Sobol' indices in practice, an algorithm that adpats to the intrinsic
randomness of the stochastic model was designed and illustrated on
simulations, where the importance of a good balance between
explorations and repetitions was observed. 
In the theoretical results, the regularization parameter, though
arbitrary, has to be fixed.  In the future, it would be of primary
interest to let it vanish as the number of explorations and
repetitions go to infinity, allowing one to get a genuine convergence
rate for the mean-squared error of the regularized estimator. A further
step would be to get a minimax low bound. Also, it would be
interesting to compare the results with multilevel Monte Carlo methods
(\cite{mycek,giles2019multilevel}). Furthermore, it could be
interesting to couple the iterative estimation approach of
\cite{gilquin} to the algorithms implemented in this study in order to
build an adaptive version which could perform estimation with respect
to a given precision.  Another estimation approach worth considering
is that based on Chatterjee's empirical correlation
coefficient~\citep{chatterjee_new_2021,gamboa_global_2022}.


\subsection*{Acknowledgment}
We thank a anonymous reviewer for their helpful comments.  We thank Clémentine
Prieur for reading the first version of this paper and for the
pertinent remarks and references she suggested for its improvement.
\renewcommand{\thesection}{\Alph{section}}
\setcounter{section}{0}

\renewcommand\thefigure{S\arabic{figure}}    
\setcounter{figure}{0}    

\newtheorem{theorema}{Theorem}[section]
\section{Proof of Proposition~\ref{thm:basic-inequality}}

%

Using convexity inequality, for all $h\in (0,1)$ and $m\ge 1$, it holds:
\begin{align*}
\E{g_h\left(\widehat{\theta}\right)-g\left(\theta\right)}^2 \leq  2\E{g_h(\widehat{\theta})-g_h\left(\mu_m\right)}^2+
2\left(g_h\left(\mu_m\right)-g\left(\theta\right)\right)^2.
\end{align*}
Applying a Taylor-Lagrange expansion to $g_h$ at points $\widehat{\theta}$ and $\mu_m$ yields:
\begin{align*}
&\E{g_h(\widehat{\theta})-g_h(\mu_m)}^2\\
&= \mathbb{E}\left(\nabla g_h\left(\mu_m\right)^{\top}\left(\widehat{\theta}-\mu_m\right)+
\frac{1}{2}\left(\widehat{\theta}-\mu_m\right)^{\top}\nabla^2g_h\left(\lambda\widehat{\theta}+(1-\lambda)\mu_m\right)\left(\widehat{\theta}-\mu_m\right)\right)^2\\
&\le 2 \mathbb{E}\left(\nabla g_h\left(\mu_m\right)^{\top}\left(\widehat{\theta}-\mu_m\right)\right)^2 +\frac{1}{2}\mathbb{E}\left(\left(\widehat{\theta}-\mu_m\right)^{\top}\nabla^2g_h\left(\lambda\widehat{\theta}+(1-\lambda)\mu_m\right)\left(\widehat{\theta}-\mu_m\right)\right)^2,
\end{align*}
for some $\lambda\in (0,1)$.  
Thus:
\begin{align*}
\frac{\E{g_h\left(\widehat{\theta}\right)-g_h\left(\mu_m\right)}^2}{\mathbb{E}\left(\nabla g_h\left(\mu_m\right)^{\top}\left(\widehat{\theta}-\mu_m\right)\right)^2}\le \quad 2+\frac{1}{2}\frac{\mathbb{E}\left(\left(\widehat{\theta}-\mu_m\right)^{\top}\nabla^2g_h\left(\lambda\widehat{\theta}+(1-\lambda)\mu_m\right)\left(\widehat{\theta}-\mu_m\right)\right)^2}{\mathbb{E}\left(\nabla g_h\left(\mu_m\right)^{\top}\left(\widehat{\theta}-\mu_m\right)\right)^2}
\end{align*}
 Let 
 \begin{equation}\label{eq:psi}
 p_{n,m}(h):=\sup_{\lambda\in [0,1]}\frac{\E{\left(\widehat{\theta}-\mu_m\right)^{\top}\nabla^2g_h\left(\lambda\widehat{\theta}+(1-\lambda)\mu_m\right)\left(\widehat{\theta}-\mu_m\right)}^2}{\E{\nabla g_h\left(\mu_m\right)^{\top}\left(\widehat{\theta}-\mu_m\right)}^2},
\end{equation}
 then 
\begin{align*}
\E{g_h(\widehat{\theta})-g_h(\mu_m)}^2\le (2+\frac{1}{2}p_{n,m}(h))V_{n,m}(h)\leq 2(1+p_{n,m}(h))V_{n,m}(h).
\end{align*}
Now, let us show that $ p_{n,m}(h)\to 0$ as $n,m\to \infty$. 
Using Cauchy-Schwarz inequality yields that: 
\begin{equation*}
 p_{n,m}(h)\leq
  \sqrt{\sup_{\lambda\in [0,1]}
    \E{
      \|\nabla^2 g_h\left(\lambda\widehat{\theta}+(1-\lambda)\mu_m\right)\|_F^4}}
\times
  \frac{\sqrt{ \E{\|\widehat{\theta}-\mu_m\|^8}}}%
  {\E{\nabla g_h\left(\mu_m\right)^{\top}\left(\widehat{\theta}-\mu_m\right)}^2}.
\end{equation*}
As $n, m\rightarrow \infty$, the first term of the product above is
bounded by a constant independent of $h$,  uniformly in $\lambda$, by Assumption~\ref{assumption_3}.  The second term is given by:
\begin{align}
  &\frac{\sqrt{ \E{\|\widehat{\theta}-\mu_m\|^8}}}
    {\E{\nabla g\left(\mu_m+h\mathbf{u}\right)^{\top}\left(\widehat{\theta}-\mu_m\right)}^2}\\
  &=\frac{\sqrt{ \E{\|\widehat{\theta}-\mu_m\|^8}}}
    {\E{\nabla g\left(\mu_m\right)^{\top}\left(\widehat{\theta}-\mu_m\right)}^2}\times \frac{\E{\nabla g\left(\mu_m\right)^{\top}\left(\widehat{\theta}-\mu_m\right)}^2}{\E{\nabla g\left(\mu_m+h\mathbf{u}\right)^{\top}\left(\widehat{\theta}-\mu_m\right)}^2}\nonumber\\
  &\leq \frac{\sqrt{ \E{\|\widehat{\theta}-\mu_m\|^8}}}
    {\E{\nabla g\left(\mu_m\right)^{\top}\left(\widehat{\theta}-\mu_m\right)}^2}\times \sup_{h\in (0,1)}\frac{\E{\nabla g\left(\mu_m\right)^{\top}\left(\widehat{\theta}-\mu_m\right)}^2}{\E{\nabla g\left(\mu_m+h\mathbf{u}\right)^{\top}\left(\widehat{\theta}-\mu_m\right)}^2}.\label{eq:ratio}
\end{align}

\begin{lemmaa}\label{lem:RatioIsOhOne}
  \begin{equation}
    \frac{\sqrt{ \E{\|\widehat{\theta}-\mu_m\|^8}}}
    {\E{\nabla g\left(\mu_m\right)^{\top}\left(\widehat{\theta}-\mu_m\right)}^2}=o(1),\label{eq:ratio-moments}
  \end{equation}
  as $n,m\to +\infty$ 
  \label{res:moments}
\end{lemmaa}

The first factor in the right-hand side in \eqref{eq:ratio} does not
depend on $h$ and is of order $o(1)$ as stated in Lemma
\ref{res:moments}.  Moreover, since $\nabla g$ is continuous in
parameter $h$, we have that the second factor in the right-hand side
of~\eqref{eq:ratio} is $O(1)$ as $n, m\to\infty$.  Therefore,
$p_{n,m}(h)=o(1)$, and hence
\begin{align*}
  \E{g_h(\widehat{\theta})-g(\theta)}^2 \leq 4\left(1+p_{n,m}(h)\right)\left(V_{n,m}(h)+B_{m}(h)^2\right),
\end{align*}
where $ B_{m}(h)=g_h(\mu)-g(\theta)$ and $\lim_{n,m\to +\infty} p_{n,m}(h)=0$.

\subsection*{Proof of Lemma ~\ref{res:moments}}

First, let us bound the numerator of the ratio in Equation \eqref{eq:ratio-moments}; we have
  \begin{equation*}
    \E{\|\widehat{\theta}-\mu_m\|^8} \le 27
    \sum_{j=1}^3\E{\|\widehat{\theta}_j-\mu_{mj}\|^8},
  \end{equation*}
  where $\widehat{\theta}_j$ and $\mu_{mj}$ denote the $j$th component of $\widehat\theta$ and $\mu_m$, respectively.
  By \cite{zygmund} and Jensen inequalities, we have for every $j$ that
  \begin{equation*}
    \E{\|\widehat{\theta}_j-\mu_{mj}\|^8} \le
    \frac{B_8}{n^4}\E{\left|\widehat\theta_{mj}^{(1)}-\mu_{mj}\right|^8},
  \end{equation*}
  where here $B_8$ is a universal constant.
  
The case $j=2$ is the simplest. 
Notice that $\mu_{m2}$ does not depend on $m$, then the expansion of $
\Esmall{|\widehat Q_m(\X^{(1)})-\mu_{m2}|^8}$ through Newton formula
yields terms of the form $\mu_{m2}^k\mathbb{E}(\widehat Q_m^{8-k}),
k=0,\dots ,8$. Using Lemma ~\ref{lemma:form} in
Section~\ref{sec:Lemma} provides that those terms are polynomial in $m^{-1}$.

Let us deal with the case $j=1$. Expanding the power $8$ through
  Newton's formula and bounding its terms yields
  \begin{equation}\label{eq:newtonBound}
    \E{\left|\widehat Q_m(\X^{(1)})^2 - \mu_{m1}\right|^8}
    \le (\mu_{m1}^8\vee 1)\binom{8}{4}\left(\E{\widehat Q_m(\X^{(1)})^{16}}+1\right).
  \end{equation}
  Denoting $\varphi(\X^{(1)},Z^{(1,k)})=Y^{(1,k)}$, we have
  \begin{align*}
    \E{\left|\widehat Q_m(\X^{(1)}) \right|^{16}}
    &= \E{\left|\frac{1}{m}\sum_{k=1}^mY^{(1,k)}
      \right|^{16}}\\
    &=  \frac{1}{m^{16}} \sum_{k_1,\dots,k_{16}=1}^m
            \E{Y^{(1,k_1)}\cdots Y^{(1,k_{16})}}.
  \end{align*}                  
  The expectation in the right-hand side is symmetric in
  $k_1,\dots,k_{16}$, and hence, from Lemma~\ref{lemma:form}, the sum
  is a polynomial in $m$ of degree $16$. Therefore, the right-hand
  side in~\eqref{eq:newtonBound} is bounded uniformly in $m$.

Let us deal with the case $j=3$. Proceeding as
  in~(\ref{eq:newtonBound}), we have
  \begin{align*}
    &\E{\left|
        \widehat{Q}_m(\X^{(1)})\widetilde{Q}_m(\widetilde{\X}^{(1)}_{\sim u},\X^{(1)}_{u})
          - \mu_{m3}\right|^8}\\
    &\le (\mu_{m3}^8\vee 1)\binom{8}{4}\left(\E{\left|
         \widehat{Q}_m(\X^{(1)})\widetilde{Q}_m(\widetilde{\X}^{(1)}_{\sim u},\X^{(1)}_{u})\right|^8}
      +1\right)\\
    &\le (\mu_{m3}^8\vee 1)\binom{8}{4}\left(\E{\frac{1}{2}
         \widehat{Q}_m(\X^{(1)})^{16}+\frac{1}{2}\widetilde{Q}_m(\widetilde{\X}^{(1)}_{\sim u},\X^{(1)}_{u})^{16}}
      +1\right),
  \end{align*}
  and this is also bounded uniformly in $m$. (Again by Lemma~\ref{lemma:form}.)

We now deal with the root of the denominator of the ratio in Equation \eqref{eq:ratio-moments}. We have
  \begin{align}
    &\mathbb{E}\left(\nabla
      g\left(\mu_m\right)^T\left(\widehat{\theta}-\mu_m\right)\right)^2\notag\\
    &=\sum_{j_1,j_2=1}^3 \nabla g(\mu_m)_{j_1} \nabla g(\mu_m)_{j_2}
                         \mathbb E(\widehat\theta-\mu_m)_{j_1}(\widehat\theta-\mu_m)_{j_2}\notag\\
    &=\frac{1}{n}\sum_{j_1,j_2=1}^3 \nabla g(\mu_m)_{j_1} \nabla g(\mu_m)_{j_2}
             \left(\mathbb E \widehat\theta^{(1)}_{mj_1}\widehat\theta^{(1)}_{mj_2}
                                  - \mu_{mj_1}\mu_{mj_2}\right).\label{eq:dslim}
  \end{align}
  The infimum of the sum in~(\ref{eq:dslim}) is reached for some $m$ and greater than zero.
Therefore, the numerator in Equation \eqref{eq:ratio-moments} is
  less than $1/n^4$ times a constant not depending on $m$ or $n$ and
  the denominator is equal to $1/n^2$ times a quantity greater than zero.
  Therefore, the supremum over $m$ of the ratio in
  Equation \eqref{eq:ratio-moments} is of order $O(n^{-2})$. The proof is
  complete.


\section{Proof of Theorem~\ref{res:main}}

  The following lemma will be needed:
  \begin{lemmaa}\label{res:inv-moments}
   For all $\alpha >0$ and all $h \in (0,1)$,
    \begin{equation}
     \lim_{n,m\to+\infty} \mathbb{E}\left[\frac{1}{
          \left(h+\widehat{\theta}_1-\widehat{\theta}_2^2\right)^{\alpha}}\right]=\frac{1}{\left(h+\var\left(Q(\X)\right)\right)^\alpha}\leq \var\left(Q(\X)\right)^{-\alpha}.
    \end{equation}
  \end{lemmaa}
  
Note that the function $g$ is infinitely differentiable over its domain $\mathcal{D}$ and then its gradient is given by:
$$
\nabla g\left(\theta_1,\theta_2,\theta_3\right)=\left(-\frac{\theta_3-\theta_2^2}{(\theta_1-\theta_2^2)^2}\,\frac{2\theta_2(\theta_3-\theta_1)}{(\theta_1-\theta_2^2)^2},\frac{1}{\theta_1-\theta_2^2}\right)^\top.
$$
Furthermore, the hessian matrix of $g$ yields:
\begin{equation*}
 \nabla^2 g(\theta_1,\theta_2,\theta_3)=\left(\begin{array}{ccc}
 \frac{2(\theta_3-\theta_2^2)}{(\theta_1-\theta_2^2)^3}& \frac{2\theta_2(\theta_1-2\theta_3+\theta_2^2)}{(\theta_1-\theta_2^2)^3} & \frac{-1}{(\theta_1-\theta_2^2)^2}  \\
 \frac{2\theta_2(\theta_1-2\theta_3+\theta_2^2)}{(\theta_1-\theta_2^2)^3}    &  \frac{2(\theta_3-\theta_1)(\theta_1+3\theta_2^2)}{(\theta_1-\theta_2^2)^3}& \frac{2\theta_2}{(\theta_1-\theta_2^2)^2} \\
   \frac{-1}{(\theta_1-\theta_2^2)^2}   & \frac{2\theta_2}{(\theta_1-\theta_2^2)^2} & 0
 \end{array}
 \right).   
\end{equation*}
For any $(\theta_1,\theta_2,\theta_3)\in\mathcal{D}$, the matrix $ \nabla^2 g(\theta_1,\theta_2,\theta_3)$ is under the form $ \nabla^2 g(\theta_1,\theta_2,\theta_3)=B(\theta_1,\theta_2,\theta_3)/(\theta_1-\theta_2^2)^3$ where $B(\theta_1,\theta_2,\theta_3)$ is the matrix:
\begin{equation*}
B(\theta_1,\theta_2,\theta_3)=\left(\begin{array}{ccc}
2(\theta_3-\theta_2^2)&2\theta_2(\theta_1-2\theta_3+\theta_2^2) & -(\theta_1-\theta_2^2) \\
 2\theta_2(\theta_1-2\theta_3+\theta_2^2)   & 2(\theta_3-\theta_1)(\theta_1+3\theta_2^2)& 2\theta_2(\theta_1-\theta_2^2) \\
 -(\theta_1-\theta_2^2)   & 2\theta_2(\theta_1-\theta_2^2)& 0
 \end{array}
 \right)   .
\end{equation*}
Notice that  $B(\theta_1,\theta_2,\theta_3)$ includes only multivariate polynomials of variables $\theta_1,\theta_2$ and $\theta_3$.

Let us check Assumption~\ref{assumption_3}. Let $\lambda\in [0,1]$ and $h\in (0,1)$. We
  have
\begin{align*}
\nabla^2 g\left(\lambda\widehat{\theta}+(1-\lambda)\mu_m+h\mathbf{u}\right) &= \frac{B\left(\lambda\widehat{\theta}+(1-\lambda)\mu_m+h\mathbf{u}\right)}{\left(h+\lambda\widehat{\theta}_1+(1-\lambda)\mu_{m1}-\left(\lambda\widehat{\theta}_2+(1-\lambda)\mu_{m2}\right)^2\right)^3}.
\end{align*}
Thus:
\newcommand{\suplambda}{\sup_\lambda}
\begin{align*}
&\E{\sup_{\lambda}\|\nabla^2 g\left(\lambda\widehat{\theta}+(1-\lambda)\mu_m+h\mathbf{u}\right)\|^4_F} \\
&=\E{\sup_{\lambda}\frac{\|B\left(\lambda\widehat{\theta}+(1-\lambda)\mu_m+h\mathbf{u}\right)\|^4_F}{\left(h+\lambda\widehat{\theta}_1+(1-\lambda)\mu_{m1}-\left(\lambda\widehat{\theta}_2+(1-\lambda)\mu_{m2}\right)^2\right)^{12}}}\\
&\leq \sqrt{\E{\sup_{\lambda}\frac{1}{\left(h+\lambda\widehat{\theta}_1+(1-\lambda)\mu_{m1}-\left(\lambda\widehat{\theta}_2+(1-\lambda)\mu_{m2}\right)^2\right)^{24}}}} \\
&\quad\times\sqrt{\E{\sup_{\lambda}\|B\left(\lambda\widehat{\theta}+(1-\lambda)\mu_m+h\mathbf{u}\right)\|^8_F}}\\ 
&\leq \sqrt{\E{\sup_{\lambda}\frac{1}{\left(h+\lambda \left(\widehat{\theta}_1-\widehat{\theta}_2^2\right)+(1-\lambda)(\mu_{m1}-\mu_{m1}^2)\right)^{24}}}} \qquad \text{ (by convexity inequality) }\\
&\quad\times\sqrt{\E{\sup_{\lambda}\|B\left(\lambda\widehat{\theta}+(1-\lambda)\mu_m+h\mathbf{u}\right)\|^8_F}}\\
&\leq \sqrt{\E{\sup_{\lambda}\frac{\lambda}{\left(h+\left(\widehat{\theta}_1-\widehat{\theta}_2^2\right)\right)^{24}}}+\sup_{\lambda}\frac{1-\lambda}{(h+\mu_{m1}-\mu_{m1}^2)^{24}}} \qquad \text{ (by convexity inequality) }\\
&\quad\times\sqrt{\E{\sup_{\lambda}\|B\left(\lambda\widehat{\theta}+(1-\lambda)\mu_m+h\mathbf{u}\right)\|^8_F}}\\
&\leq \sqrt{\E{\frac{1}{\left(h+\left(\widehat{\theta}_1-\widehat{\theta}_2^2\right)\right)^{24}}}+\frac{1}{(h+\mu_{m1}-\mu_{m1}^2)^{24}}} \\
&\quad\times\sqrt{\E{\sup_{\lambda}\|B\left(\lambda\widehat{\theta}+(1-\lambda)\mu_m+h\mathbf{u}\right)\|^8_F}}\\
&\leq \sqrt{\E{\frac{1}{\left(h+\left(\widehat{\theta}_1-\widehat{\theta}_2^2\right)\right)^{24}}}+\sup_{h\in (0,1)}\frac{1}{(h+\mu_{m1}-\mu_{m1}^2)^{24}}} \\
&\quad\times\sqrt{\sup_{h\in (0,1)}\E{\sup_{\lambda}\|B\left(\lambda\widehat{\theta}+(1-\lambda)\mu_m+h\mathbf{u}\right)\|^8_F}}.
\end{align*}
One should remark that $\sup_{h\in  (0,1)}\frac{1}{(h+\mu_{m1}-\mu_{m1}^2)^{24}}\leq \frac{1}{(\V{\Ec{\varphi(\X,Z)}{\X}})^{24}} <+\infty$.  Moreover,  the matrix $B$ is composed with polynomials of three variables. Since $\E{Q(\X)^{16}}<+\infty$ then by using Lemma~\ref{lemma:form} and by continuity of polynomial functions, it yields that \\
$\sup_{h\in (0,1)}\E{\sup_{\lambda\in [0,1]}\|B\left(\lambda\widehat{\theta}+(1-\lambda)\mu_m+h\mathbf{u}\right)\|^8_F}$ is bounded.  Finally, 
by relying on Lemma \ref{res:inv-moments},
$\mathbb{E}(h+\widehat{\theta}_1-\widehat{\theta}_2^2)^{-24}$ is a
bounded by $\frac{1}{(\V{\Ec{\varphi(\X,Z)}{\X}})^{24}}$ as
$n,m\rightarrow +\infty$.  Therefore,   Assumption \ref{assumption_3}
is satisfied.

\subsection*{Proof of Lemma~\ref{res:inv-moments}}
Let $h\in (0,1)$ be fixed.  
The function $\beta_h: x\mapsto 1/(h+x)^\alpha$ is continuously differentiable such that its first derivative is uniformly bounded on $\mathbb{R}_+$ by $1/h^{\alpha+1}$ then it is Lipschitz. Therefore:
\begin{align*}
\mathbb{E}\left(\beta_h(\widehat{\theta}_1-\widehat{\theta}_2^2)-\beta_h(\theta_1-\theta_2^2)\right)^2 &\le \frac{1}{h^{2\alpha+2}}\mathbb{E}\left(\widehat{\theta}_1-\widehat{\theta}_2^2-\theta_1+\theta_2^2\right)^2\\
&\le \frac{2}{h^{2\alpha+2}}\left(\var(\widehat{\theta}_1)+\var(\widehat{\theta}_2^2)+\left(\mathbb{E}(\widehat{\theta}_2^2)-\theta_2^2\right)^2\right),
\end{align*}
using convexity inequality.  Based  on Marcinkiewicz-Zygmund inequality (see Theorem \ref{thm2}) and Lemma \ref{lemma:form}, it follows that
$$
\lim_{n,m\to\infty}\left((\var(\widehat{\theta}_1)+\var(\widehat{\theta}_2^2)+\left(\mathbb{E}(\widehat{\theta}_2^2)-\theta_2^2\right)^2\right)=0.
$$
Straightforwardly:
$$
\lim_{n,m\to \infty}\mathbb{E}\left(\frac{1}{(h+(\widehat{\theta}_1-\widehat{\theta}_2^2)^\alpha}\right)=\frac{1}{(h+\theta_1-\theta_2^2)^\alpha}\leq \frac{1}{\var(Q(\X))^\alpha}.
$$

\section{A lemma}\label{sec:Lemma}

\begin{lemmaa}
\label{lemma:form}
Let $\X^{(1)}, \cdots ,\X^{(n)}$ be $n$ i.i.d. copies of $\X$ and
$Z^{(1,1)}, \cdots ,Z^{(n,m)}$ be $n\times m$ i.i.d.  copies
of $Z$ such that $(\X^{(1)}, \cdots ,\X^{(n)})$ and $(Z^{(1,1)},
\cdots ,Z^{(n,m)})$ are independent.  Then, for all $q\in \mathbb{N}$:
 $m^{-q}\E{\sum_{k=1}^m\varphi(\X^{(1)},Z^{(1,k)})}^q$ is polynomial in $m^{-1}$ of degree $q-1$ with constant $\E{\Ec{\varphi(\X,Z)}{\X}}^q$.

\end{lemmaa}

 It holds that:
\begin{align*}
\E{\frac{1}{m}\sum_{k=1}^m\varphi(\X^{(1)},Z^{(1,k)})}^q &= \frac{1}{m^q}\E{\sum_{k_1=1}^m\cdots \sum_{k_q=1}^m\varphi(\X^{(1)},Z^{(1,k_1)})\cdots \varphi(\X^{(1)},Z^{(1,k_q)}) }\\
&=\frac{1}{m^q}\sum_{k_1=1}^m\cdots \sum_{k_q=1}^m\E{\varphi(\X^{(1)},Z^{(1,k_1)})\cdots \varphi(\X^{(1)},Z^{(1,k_q)}) }.
\end{align*}


  Denote by $\lambda:\{1,\dots,m\}^q\to\mathbb{N}$ the map which with
  each $\mathbf k:=(k_1,\dots,k_q)$ associates the number of distinct indices among
  $k_1,\dots,k_q$. If $1\leq l\leq q$ then denote by  $\rho_l:\lambda^{-1}(l)\to\{1,\dots,q\}^l$ the map which with each $\mathbf k\in\lambda^{-1}(l)$ associates $(r_1,\dots,r_l)$, 
 where $r_i=|\{j:k_j=k_{j_i}\}|$ for every
$i=1,\dots,l$ and $k_{j_1},\dots,k_{j_l}$ are the distinct indices
found among $k_1,\dots,k_q$. Obviously, $r_1+\cdots+r_l=q$. We have
\begin{align}\label{eq:newsum}
  &\sum_{k_1=1}^m\cdots \sum_{k_q=1}^m\E{\varphi(\X^{(1)},Z^{(1,k_1)})\cdots \varphi(\X^{(1)},Z^{(1,k_q)}) }\notag\\
  &=\sum_{k_1=1}^m\cdots \sum_{k_q=1}^mf(\mathbf k)\notag\\
  &=\sum_{l=1}^q
    \left(\sum_{(r_1,\dots,r_l)\in \{1,\dots,q\}^l : r_1+\cdots+r_l=q}
    \left(\sum_{\mathbf k\in\lambda^{-1}(l):\rho_l(\mathbf k)=(r_1,\dots,r_l)}
    f(\mathbf k)\right)\right).
\end{align}
Now, since
\begin{align*}
  &\E{\varphi(\X^{(1)},Z^{(1,k_1)})\cdots \varphi(\X^{(1)},Z^{(1,k_q)})}\\
   &= \E{\varphi(\X^{(1)},Z^{(1,k_{j_1})})^{r_1}\cdots
       \varphi(\X^{(1)},Z^{(1,k_{j_l})})^{r_l}}\\
   &=\E{\prod_{s=1}^l\Ec{\varphi(\X^{(1)},Z^{(1,k_{j_s})})^{r_s}}{\X^{(1)}}}
\end{align*}
is symmetric in $r_1,\dots,r_l$, it holds that
\begin{align*}
  &\sum_{(r_1,\dots,r_l)\in \{1,\dots,q\}^l : r_1+\cdots+r_l=q}\left(
    \sum_{\mathbf k\in\lambda^{-1}(l):\rho_l(\mathbf k)=(r_1,\dots,r_l)}
    f(\mathbf k)\right)\\ 
  &= c(l,(r_1,\dots,r_l),m) 
    \E{\prod_{s=1}^l\Ec{\varphi(\X^{(1)},Z^{(1,k_{j_s})})^{r_s}}{\X^{(1)}}}
\end{align*}
where 
\begin{multline}\label{eq:coef}
  c(l,(r_1,\dots,r_l),m) =
  \binom{q}{r_1}\binom{q-r_1}{r_2}\cdots\binom{q-r_1-\cdots-r_{l-1}}{r_l}\\
    m(m-1)\cdots(m-l+1).
\end{multline}
Notice
that the expression in the right-hand side of~\eqref{eq:coef} is
invariant by permutation of $r_1,\dots,r_l$. 
Therefore, the sum~\eqref{eq:newsum} is a polynomial in $m$ of degree $q$ with constant zero and hence $\E{\frac{1}{m}\sum_{k=1}^m\varphi(\X^{(1)},Z^{(1,k})}^q$ is
a polynomial in $\frac{1}{m}$ of degree $q-1$ with constant $\lim_{m\rightarrow +\infty}\E{\frac{1}{m}\sum_{k=1}^m\varphi(\X^{(1)},Z^{(1,k})}^q=\E{\Ec{\varphi(\X,Z)}{\X}}^q$.

\section{The Marcinkiewicz-Zygmund inequality}
\begin{theorema}[\cite{zygmund}]
Let $U_1,\cdots ,U_n$ be i.i.d. random variables such that $\mathbb{E}(U_1)=0$ and $\mathbb{E}\vert U_1\vert^q<+\infty$, where $1\leq q<+\infty$. There exist $A_q$ and $B_q$ depending only on $q$ such that:
\begin{equation*}
A_q\mathbb{E}\left(\left(\sum_{i=1}^n\left|U_i\right|^2\right)^{\frac{q}{2}}\right)\leq \mathbb{E}\left(\left|\sum_{i=1}^nU_i\right|^q\right)\leq B_q\mathbb{E}\left(\left(\sum_{i=1}^n\left|U_i\right|^2\right)^{\frac{q}{2}}\right)\label{thm2}
\end{equation*}
Furthermore, there exists $C_q$ independent from $n$ such that:
\begin{equation}
\mathbb{E}\left(\left|\frac{1}{n}\sum_{i=1}^nU_i\right|^q\right)\leq \frac{C_q}{n^{\frac{q}{2}}}. \label{moment_rate}
\end{equation}
\end{theorema}


\bibliographystyle{plainnat}
\bibliography{mabiblio.bib}

\end{document}